\documentclass[english]{article}

\usepackage{amsmath,amsfonts,amsthm,amssymb}
\usepackage{graphicx}
\usepackage{soul}
\usepackage{epstopdf,ulem}
\usepackage{enumerate}
\usepackage{booktabs}
\usepackage{hyperref}
\usepackage{fancyhdr}
\usepackage{tgschola}
\usepackage[overload]{empheq} 
\usepackage{color}
\usepackage{float}
\usepackage[all]{xy}
\usepackage{tikz}
\usetikzlibrary{arrows}
\definecolor{mosco}{cmyk}{1,0,0,0}
\usepackage{graphicx}
\usepackage{pstricks}
\usepackage{pstricks-add}
\usepackage{pstcol,pst-3d,pst-char,pst-coil,pst-eps,pst-fill,pst-grad,pst-node,pst-plot,pst-text,pst-tree,pst-math}

\newcommand{\N}{\mathcal{N}}
\newcommand{\Xb}{\mathbf{X}}

\newcommand{\Gb}{\mathbf{G}}
\newcommand{\Hb}{\mathbf{H}}
\newcommand{\Eb}{\mathbf{E}}
\newcommand{\Phib}{\mathbf{\Phi}}

\newtheorem{theorem}{Theorem}
\newtheorem{proposition}{Proposition}

\newtheorem{remark}{Remark}

\title{Pheromone trapping for control of Asian citrus psyllid, \textit{Diaphorina citri} Kuwayama (Hemiptera: Liviidae)}
\vspace{1cm}

\author{ Daiver Cardona-Salgado$^{1}$, \ Yves Dumont$^{2,3,4}$, \ Olga Vasilieva$^{5}$\footnote{Corresponding author: olga.vasilieva@correounivalle.edu.co} \vspace{5mm} \\
$^1$\small Department of Mathematics, Universidad Autonoma de Occidente, Cali, Colombia \\
$^2$ \small CIRAD, UMR AMAP, P\^{o}le de Protection des Plantes, F-97410 St Pierre, R\'{e}union island, France.\ \\ $^3$ \small AMAP, University of Montpellier, CIRAD, CNRS, INRAE, IRD, Montpellier, France \\
$^4$ \small Department of Mathematics and Applied Mathematics, University of Pretoria, Pretoria, South Africa \\
$^5$\small Department of Mathematics, Universidad del Valle, Cali, Colombia \\
}

\date{\today}

\begin{document}
\maketitle
\begin{abstract}
We study the impact of pheromone control against the Asian citrus psyllid, \textit{Diaphorina citri}, a principal vector of diseases in citrus cultures. The model is expressed as a piecewise smooth ODE system, and its long-term behavior is analyzed. In particular, through qualitative analysis and applying an open-loop control approach, we identify the threshold in terms of two external parameters related to the pheromone traps, the amount of pheromones to be released and the male-killing rate, to ensure local elimination of the wild psyllid population. We also show that a feedback control with periodic assessments of the wild population sizes is applicable and then deduce that a mixed-type control, combining the open- and closed-loop control approaches, provides the best results. We present several simulations to illustrate our theoretical findings and to estimate the minimal amount of pheromones and time needed to reach the local elimination of wild psyllids. Finally, we discuss possible implementations  of our results as a part of Integrated Pest Management programs. \\ \\

\end{abstract}

%
\baselineskip 7mm

\section{Introduction}
\label{sec-intro}

The Asian citrus psyllid (ACP) \cite{Aidoo2022,Grafton-Cardwell2013}, \textit{Diaphorina citri} Kuwayama, is the most important pest of citrus cultures because it is the main vector of \textit{Candidatus} Liberibacter spp., the bacterium that cause huanglongbing (HLB), the citrus greening disease \cite{Bove2006}, impacting several places around the world, and, in particular, Colombia \cite{Angel2014} and also La R\'{e}union, a French overseas department. When uninfected psyllids feed on an infected citrus tree, they acquire the bacterium. Subsequently, when they feed on healthy trees, they can transmit the bacterium, thereby spreading the disease. Note also that in La R\'eunion another psyllid, the African citrus psyllid, \textit{Trioza erytreae}, has been identified as an efficient vector of \textit{Candidatus} Liberibacter asiaticus \cite{Reynaud2022}. Since there is no cure for infected trees, several control strategies have been developed, including the removal and destruction of infected trees to prevent further spread, the use of insecticides to control the Asian citrus psyllid population, quarantine measures to limit the movement of infected plant material, research into disease-resistant citrus varieties through breeding programs, development of early detection methods to identify infected trees, etc. So far, only in La R\'eunion the biological HBL control was achieved successfully in the late 1970s \cite{Aubert1984}. Indeed, in La R\'eunion, the two psyllid vectors have been controlled with hymenopteran psyllid parasites: \textit{Tamarixia radiata} introduced from India against \textit{D. citri} and \textit{Tamarixia dryi}, from South Africa, against \textit{Trioza erytreae} \cite{Aubert1984}. \textit{D. citri}  started to be reported in the Caribbean basin in the late 1990s \cite{Halbert2004}, before being first officially reported in Colombia in 2007 \cite{Ebratt2011,Ram2018}. ACP Biological control methods started in Colombia, using, for instance, natural enemies of \textit{D. citri}, like \textit{T. radiata} \cite{Ebratt2011}, and others collected in the department of Valle del Cauca, Colombia \cite{Kondo2015}. In Colombia, single insecticides and insecticide rotations have also been tested against \textit{D. citri} \cite{Ram2018}. Meanwhile, in Brazil, ACP biological control with sex pheromones \cite{Zanardi2018} is under study.

Mathematical modeling is now a common tool to study (biological) control strategies against pests \cite{Tapi2020} and vectors \cite{Anguelov2017}. In particular, several models have been developed and studied to control the spreading of HLB: see \cite{Taylor2016} for an overview and references therein. The majority of these models are epidemiological models based on vector-borne disease models developed for mosquitoes. In \cite{Gao2021}, the authors developed an ACP continuous population model to study the effect of physiological and behavioral resistance and investigate the existence of threshold conditions for extinction. Discrete ACP models for each stage (eggs, nymphs, and adults) have been developed in \cite{Milosavljevic2018} to study the impact of environmental parameters, habitat, and natural enemies on the ACP dynamics in an urban area in California. However, these phenological models are degree-day models, i.e., based on a temperature accumulation, and thus well adapted to study population accumulation of \textit{D. citri} and the effect of temperature. In this paper, we consider a piecewise smooth modeling approach to study the impact of sex pheromone control.

The outline of the paper is as follows. In Section \ref{sec-model}, we propose a sex-structured mathematical model that encompasses only the population of adult Asian citrus psyllids. The model is formulated as a piecewise smooth dynamical system in continuous time. In Section \ref{sec-analysis}, dedicated to the qualitative analysis of the proposed model, the long-term evolution of the natural ACP population dynamics is studied, and the underlying stability properties of the piecewise smooth dynamical system are established. The proposed model is further amended in Section \ref{sec-discussion} with external control actions of pheromone traps: attraction and direct removal of male insects that induce mating disruption targeting to reduce the future offspring. These intervention measures are modeled by two external parameters, the male-killing rate  and the strength of lure. The choice of these two parameters may result in two outcomes: the suppression or elimination of the local ACP population. To reach one of these goals, the open-loop and closed-loop operational control modes are suggested and validated in Subsections \ref{subsec-open} and \ref{subsec-closed}, respectively. Section \ref{sec-num} provides numerical simulation illustrating the open-loop and closed-loop control approaches. Finally, Section \ref{sec-concl} summarizes the main results of our work.

\section{Natural population dynamics of \textit{Diaphorina citri}}
\label{sec-model}

Asian citrus psyllids are small (2.7 to 3.3 mm long) jumping and flying insects that live on citrus trees and feed on young stems, sprouts, and leaves during all stages of development. The psyllid's life cycle includes an immature phase (consisting of the egg stage and five nymphal instars) followed by the adult stage (imago) of sexually matured insects, males or females. Oviposition and development of immature \textit{D. citri} elapse on young, tender flush leaves where the nymphs remain almost docile while feeding on the tissue of young leaves and stems until turning into adults \cite{Hall2020}.

In this section, we propose a sex-structured mathematical model that encompasses only the population of adult Asian citrus psyllids (ACP), \textit{Diaphorina citri}, even though the ACP life cycle also includes the immature phase (consisting of eggs and five nymphal instars). The model is based on this insect species' behavioral and biological features, and particular attention is paid to the ACP mating behavior. On the other hand, oviposition can also be reduced by the continuous presence of males seeking matings since this particular species (\textit{D. citri}) exhibits a male-biased operational sex ratio \cite{Lubanga2018}, meaning that there are more sexually active males than sexually receptive females.

Laboratory and field observations show that fertilized female psyllids become temporarily unavailable for mating and try to avoid males when they are ready for oviposition \cite{Mankin2020}. After oviposition, such female insects again exhibit receptiveness for mating. Thus, the female psyllids usually mate intermittently during their lives to keep an adequate amount of viable sperm and be able to lay eggs throughout their lives whenever young leaves and stems are present.

To mimic this mating behavior, we divide the total population of adult psyllids into three disjoint compartments or population classes, namely:
\begin{itemize}
\item	
$M(t)$  --  the number or density of male insects at the moment $t$.
\item	
$A(t)$  -- the number or density of female insects available for mating at the moment $t$.
\item
$U(t) $  -- the number or density of fertilized female insects at the moment $t$ (they avoid mating while preparing for oviposition).
\end{itemize}

\begin{figure}[t]
\begin{center}
	\pagestyle{empty}
	\thispagestyle{empty}
	\hspace{-2.5cm}
  \centering
	\xymatrix{
	*+<1cm>[F]{A}\ar@{->}[dr] \ar[d]_-{\delta} & & *+<1cm>[F]{M}\ar@{-->}[dl]  \ar[d]_-{\mu}\\
	&{\nu \min \left\{ \frac{\gamma M}{A}, 1  \right\}}\ar[d] & \\
	& *+<1cm>[F]{U}\ar@{->}[uul]^{\eta}  \ar[d]_-{\delta} \ar@[mosco]@/^4cm/[luu]+L_{(1-r) \rho e^{-\sigma (M+A+U)}} \ar[d]_-{\delta} \ar@[mosco]@/^-4cm/[uur]+R^{r \rho e^{-\sigma (M+A+U)}} & \\
	& & \\
	}
\end{center}
\caption{Flow diagram of the natural dynamics of \textit{Diaphorina citri} described by system \eqref{sys} \label{fig-1}}
\end{figure}
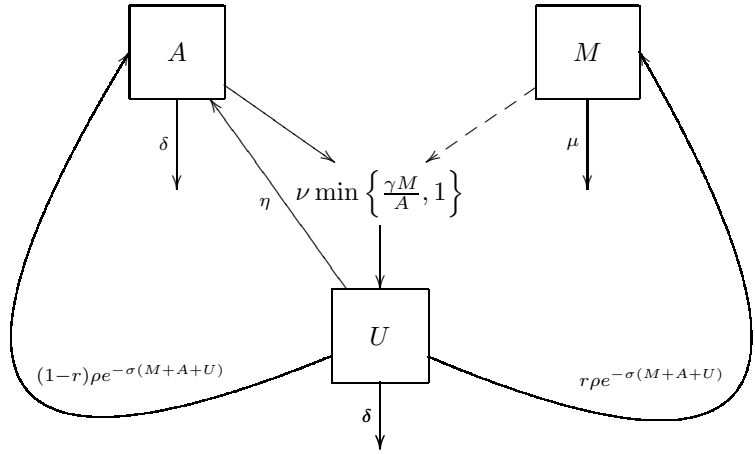

Thus, $F(t):=A(t)+U(t)$ constitutes the total population of female psyllids. It is also supposed that all male insects $M(t)$ are available for mating anytime and remain sexually active during their lifetime.

Following the approach of \cite{Anguelov2017,Tapi2020}, and according to the flow diagram provided in Figure \ref{fig-1}, we derive the following ODE system to describe the population dynamics of adult psyllids.

\begin{subequations}
\label{sys}
\begin{align}[left = \empheqlbrace\,]
\label{sys-M}
\frac{d M}{dt}& = r \rho U e^{-\sigma (M+A+U)} - \mu M \\
\label{sys-A}
\frac{d A}{dt}& = (1-r) \rho U e^{-\sigma (M+A+U)} - \nu \min \left\{ \frac{\gamma M}{A}, 1  \right\} A +  \eta U - \delta A \\
\label{sys-U}
\frac{d U}{dt}& = \nu \min \left\{ \frac{\gamma M}{A}, 1  \right\} A - \eta U - \delta  U
\end{align}
\end{subequations}

\noindent
with nonnegative initial conditions

\begin{equation}
\label{icon}
 M(0)= M_0, \quad A(0)=A_0, \quad U(0)=U_0.
 \end{equation}

The constant parameters included in the model \eqref{sys} are all positive, and their concise definitions, as well as numerical values in simulations (Section \ref{sec-num}), are summarized in Table \ref{tab-1}. Notably, these parameter values are borrowed from a field study on Valencia sweet orange tree (\textit{Citrus sinensis}) with Rangpur lime as a rootstock (\textit{Citrus limonia}).

In the equations \eqref{sys-M} and \eqref{sys-A}, we denote by $r$ and $(1-r)$ with $r \in (0,1)$ the proportion of male and female psyllids emerging from the immature stage and entering the compartments of males and receptive females, respectively.

\begin{table}[t]
\centering
\caption{Parameters of the model \eqref{sys} with some values corresponding to the field study performed on Valencia sweet orange tree (\textit{Citrus sinensis}) with Rangpur lime (\textit{Citrus limonia}) taken as a rootstock \label{tab-1}}
{\small
\begin{tabular}{clccc}
\toprule
\textbf{Parameter}	& \textbf{Description}	& \textbf{Value} & \textbf{Unit} & \textbf{References} \\
\midrule
$r$      & primary sex ratio                                  & $0.41$    & --                              & \cite{Perez2017} \\
$\rho$   & mean no. of eggs produced by one female per day    & $6.352$   & day$^{-1}$ & \cite{Perez2017} \\
$\sigma$ & characteristic of eggs survival to the adult stage & $0.001$    & individual$^{-1}$ & assumed \\
$\mu$    & natural mortality rate for males                   & $0.021$  & day$^{-1}$ & \cite{Perez2017} \\
$\delta$ & natural mortality rate for females                 & $0.023$   & day$^{-1}$ & \cite{Perez2017} \\
$\gamma$ & females fertilized by a single male                & $1.2$     & --       & \cite{Perez2017} \\
$\nu$    & transfer rate from $A$ to $U$                      & $1/4$       & day$^{-1}$ & estimated from \cite{Zanardi2018} \\
$\eta$   & transfer rate from $U$ to $A$                      & $1$  & day$^{-1}$ & - \\
\bottomrule
\end{tabular}
}
\end{table}

The parameter $\rho>0$ stands for the mean number of eggs produced on average per day by one female psyllid from the class $U$. At the same time, the exponential factor in the recruitment terms of equations \eqref{sys-M} and \eqref{sys-A} expresses the eggs' survival to adulthood while they pass through five nymphal instars. The parameter $\sigma >0$ in the exponential factor may be seen as the ratio $\sigma = \beta/ K$ between $\beta$, a quantity characterizing the transition of immature insects into adults under density dependence and nymphal competition for food resources, and a carrying capacity $K$. The latter is typically proportional to the capacity of available breeding sites (young stems, sprouts, and leaves) that also provide food for all nymphal stages and adults (males and two classes of females).

Natural mortality rates for males and females ($F=A+U$) are denoted by $\mu$ and $\delta$, respectively, and correspond to the inverses of their average lifespans ($1/\mu$ and $1/\delta$ days, respectively). Some studies report that female psyllids live longer than males (see \cite{Hall2020,Perez2017} and more detailed references therein), so we suppose in the sequel that $\mu \geq \delta$.

Further, we assume that a receptive female $A$ needs to mate once or more to pass into the class $U$ of eggs-laying females and be able to reproduce. The conversion of mating females $A$ into eggs-laying females $U$ is modeled by the mating term $\nu \min \Big\{ (\gamma M)/A, 1  \Big\} A$ that appears in equations \eqref{sys-A} and \eqref{sys-U}. In this term, $\gamma \geq 1$ expresses the mean number of females a single male can fertilize. Furthermore, the parameter $\nu$ can be viewed as the effective mating (or contact) rate that results in successful fertilization of the female leading to her readiness for oviposition. In other words, it is assumed that a sexually mature female becomes ready for oviposition after $1/\nu$ days from exhibiting receptiveness and completing at least one mating. The latter is valid only if there are enough males so all females from class $A$ can mate at least once. However, if male psyllids are scarce, then only a proportion $(\gamma M)/A$ of mate-seeking females $A$ can get fertilized and pass into the eggs-laying class $U$ for further reproduction. Moreover, after completing the oviposition, a female psyllid $U$ becomes receptive to mating again after $1/\eta$ days and moves back to $A$-class.

Notably, the three-dimensional model \eqref{sys} has been designed by merging two modeling approaches. Namely, we have used as a basis the two-dimensional sex-structured model initially developed by Bliman \textit{et al.} \cite{Bliman2019} for mimicking the population dynamics of any pest or disease vector population. It is also well known that female psyllids must re-mate after each oviposition to enhance fertility and continue laying viable eggs. To mimic this process, we have introduced a separate female class $A$ gathering all female psyllids available for mating, like in \cite{Anguelov2017}. To model the  re-mating process, we have employed the so-called ``mating function'' ($ \min \big\{ (\gamma M)/A, 1  \big\}$), which was proposed initially by Barclay \& van den Driessche \cite{Barclay1983} for discrete-time models, and further adapted to continuous-time models by Anguelov \textit{et al.} \cite{Anguelov2017} (see also \cite{Tapi2020}).

Using the approach developed in \cite{Anguelov2017}, system \eqref{sys} can be written in the  form

\begin{equation}
\label{sys-pws}
\frac{d \Xb}{dt} = \Phib(\Xb) := \left\{ \begin{array}{ccc} \Phib_1(\Xb) & \text{if} & \gamma M \geq A \\[2mm] \Phib_2(\Xb) & \text{if} & \gamma M \leq A \end{array} \right.,
\end{equation}

\noindent
where $\Xb:= (M,A,U) \in \mathbb{R}^3_+$ and

\begin{align}
\label{F1}
\Phib_1 (\Xb) & = \begin{pmatrix} r \rho U e^{-\sigma (M+A+U)} - \mu M \\[3mm]
(1-r) \rho U e^{-\sigma (M+A+U)} - \nu A +  \eta U - \delta A \\[2mm]
\nu A - \eta U - \delta  U
\end{pmatrix}, \\ & \notag \\
\label{F2}
\Phib_2 (\Xb) & = \begin{pmatrix} r \rho U e^{-\sigma (M+A+U)} - \mu M \\[3mm]
(1-r) \rho U e^{-\sigma (M+A+U)} - \nu \gamma M +  \eta U - \delta A \\[2mm]
\nu \gamma M - \eta U - \delta  U \end{pmatrix}.
\end{align}

\begin{figure}[t]
\centering
\includegraphics[width=9cm]{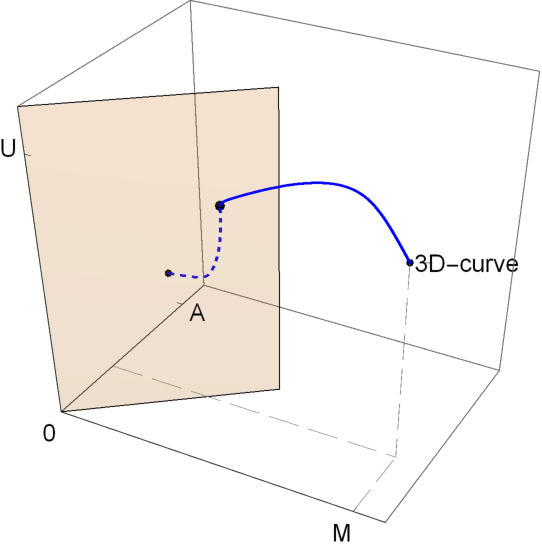}
\caption{A piecewise smooth solution $\big(M(t),A(t),U(t) \big)$ to the system \eqref{sys} drawn as a parametric 3D-curve that crosses the switching plane $\mathcal{P}_s$ (shadowed area)}
 \label{fig-3}
\end{figure}

Following definitions given in \cite{DiBernardo2008}, the dynamical system defined by \eqref{sys}, \eqref{sys-pws} can be considered as a \textit{piecewise smooth (PWS) system} with the switching manifold defined by the plane

\[ \mathcal{P}_s := \big\{ (M,A,U) \in \mathbb{R}^3_+: \ \gamma M = A \big\} \]
because any point $\tilde{\Xb}=(M,A,U) \in \mathcal{P}_s$ satisfies the relationship $\Phib_1 \big(\tilde{\Xb} \big) = \Phib_2 \big(\tilde{\Xb} \big)$. Even though the first derivatives of $\Phib$ in \eqref{sys-pws} have a jump discontinuity across the switching plane $\mathcal{P}_s$, their one-side limits are finite, and the jumps are bounded. Therefore, the overall vector field $\Phib$ is continuous and piecewise smooth for all $\Xb \in \mathbb{R}^3_+$, meaning that the right-hand side of the dynamical system \eqref{sys}, \eqref{sys-pws} is Lipschitz. The latter guarantees the existence and uniqueness of a piecewise smooth solution to the initial-value problem \eqref{sys}-\eqref{icon}. Figure \ref{fig-3} gives an example of the piecewise smooth solution $\big(M(t),A(t),U(t) \big)$ to the system \eqref{sys} in the form of a parametric 3D-curve that crosses the switching plane $\mathcal{P}_s$.

Let us denote by $\Xb(t;\Xb_0)$ the solution of \eqref{sys} engendered by the initial condition $\Xb_0:=\big( M_0,A_0,U_0 \big).$ If $\Xb_0 \in \mathbb{R}^3_+$ then it is easy to show that $\Xb(t;\Xb_0) \in \mathbb{R}^3_+$. In effect, it is fulfilled that

\[ \left. \frac{d M}{dt} \right|_{M=0}  \geq 0, \qquad  \left. \frac{d A}{dt} \right|_{A=0} \geq 0, \qquad  \left. \frac{d U}{dt} \right|_{U=0} \geq 0. \]
Therefore, the positive invariance of $\mathbb{R}^3_+$ becomes obvious and we have $\Xb(t;\Xb_0) \geq 0$ for all $t \geq 0$ whenever $\Xb_0 \in \mathbb{R}^3_+$.

Furthermore, we can establish the following result related to the uniform ultimate bound of all solutions to the PWS system \eqref{sys}.
\begin{proposition}
\label{prop1}
There exists a compact absorbing set $\Omega \subset \mathbb{R}^3_+$ that attracts all the solutions of the PWS system \eqref{sys} engendered by any initial condition $\big( M_0,A_0,U_0 \big) \in \mathbb{R}^3_+$.
\end{proposition}
\begin{proof}
First, we note that along the trajectories of \eqref{sys}, it is fulfilled that

\[ \frac{d(M + A + U)}{dt} = \rho U e^{-\sigma (M + A + U)} - (\mu M + \delta A + \delta U) \leq (M + A + U) \Big[ \rho e^{-\sigma (M + A + U)} - \min \{ \mu, \delta \} \Big]
\]
Therefore,

\[ M(t) + A(t) + U(t) \leq \max \Big\{ M_0 + A_0 + U_0, \widehat{P} \Big\} \]
where

\[ \widehat{P}:= \frac{1}{\sigma} \ln \left( \frac{\rho}{\min \{ \mu, \delta \}} \right) \]
stands for the carrying capacity of the Ricker differential equation $P'(t)=P(t) \Big[ \rho e^{-\sigma P(t)} - \min \{ \mu, \delta \} \Big].$ Thus, the compact set

\[ \Omega := \Big\{ (M,A,U) \in \mathbb{R}^3_+: \ 0 \leq M + A + U \leq \widehat{P} \Big\} \]
is invariant in the sense that any solution of \eqref{sys} engendered by $\big( M_0,A_0,U_0 \big) \in \Omega$ remains in $\Omega$ for all $t \geq 0.$ Moreover, $\Omega$ attracts all the trajectories engendered by $\big( M_0,A_0,U_0 \big) \in \mathbb{R}^3_+ \setminus \Omega$ and there is a finite time $\hat{t} >0$ such that $\big(M(\hat{t}), A(\hat{t}), U(\hat{t}) \big) \in \Omega$ for all $t \geq \hat{t}.$ In other words, $\Omega$ constitutes the absorbing set of the PWS system \eqref{sys}, and its trajectories engendered by any initial condition $\big( M_0,A_0,U_0 \big) \in \mathbb{R}^3_+$ are uniformly ultimately bounded.
\end{proof}

Once the well-posedness of the PWS system \eqref{sys} is formally established, we proceed to study its stability by applying the methodology employed in \cite{Anguelov2017,Tapi2020}.

\section{Qualitative analysis of the PWS system \eqref{sys}}
\label{sec-analysis}

The switching plane $\mathcal{P}_s$ divides the positive octant $\mathbb{R}^3_+$ into two disjoint regions:
\begin{enumerate}
\item
The \textit{male abundance region}

\[ \mathbb{M}_a := \Big\{ (M,A,U) \in \mathbb{R}^3_+: \  \gamma M > A \Big\}, \]
where the vector field $\Phib_1(\Xb)$ defined by \eqref{F1} takes action, that is,

\begin{equation}
\label{sys-F1}
\frac{d \Xb}{d t} = \Phib_1(\Xb).
\end{equation}
\item
The \textit{male scarcity region}

\[ \mathbb{M}_s := \Big\{ (M,A,U) \in \mathbb{R}^3_+: \  \gamma M < A \Big\}, \]
where the vector field $\Phib_2(\Xb)$ defined by \eqref{F2} takes action, that is,

\begin{equation}
\label{sys-F2}
\frac{d \Xb}{d t} = \Phib_2(\Xb).
\end{equation}
\end{enumerate}
As shown in Figure \ref{fig-3}, $\mathbb{M}_a$-region is in front of the switching plane $\mathcal{P}_s$ (shadowed area), whereas $\mathbb{M}_s$-region is behind $\mathcal{P}_s$. Both systems \eqref{sys-F1} and \eqref{sys-F2} have smooth right-hand sides. Their stability properties can be studied separately, at least to understand better the overall dynamics of the original PWS system \eqref{sys} whose long-term behavior is richer and more complex than that of the two smooth ODE systems \eqref{sys-F1} and \eqref{sys-F2} when considered separately.

Generally speaking, a solution $\Xb \big( t;\Xb_0 \big)$ with $\Xb_0 \in \mathbb{M}_s$ may remain in $\mathbb{M}_s$ or may enter the region $\mathbb{M}_a$ by crossing the switching plane $\mathcal{P}_s$ and then remain there (this situation is illustrated in Figure \ref{fig-3}). It is also not excluded that the mentioned solution leaves the region $\mathbb{M}_a$ and then returns or moves cyclically across the plane $\mathbb{P}_s$ (periodic or chaotic behavior). Similar behavior options also apply to solutions engendered by $\Xb_0 \in \mathbb{M}_a$. Therefore, for analyzing the long-term behavior of solutions of the PWS system \eqref{sys}, \eqref{sys-pws}, the first step will be to identify the equilibria of \eqref{sys-F1} and \eqref{sys-F2} and then to study their stability properties separately. Further, we will explore their possible connections and relations with the equilibria of the PWS system \eqref{sys}, \eqref{sys-pws}.

In this context, it is useful to recall some definitions related to the classification of equilibria a PWS system may possess \cite{DiBernardo2008}.
On the one hand, a point $\Xb^{*} \in \mathbb{R}^3_+$ satisfying either

\[ \Phib_1 \big( \Xb^{*} \big) = \mathbf{0} \quad \text{and} \quad \Xb^{*} \in \mathbb{M}_a \]
or

\[ \Phib_2 \big( \Xb^{*} \big) = \mathbf{0} \quad \text{and} \quad \Xb^{*} \in \mathbb{M}_s \]
is referred to as a \textit{regular equilibrium} of the PWS system \eqref{sys}, \eqref{sys-pws}. On the other hand, a point $\Xb^{*} \in \mathbb{R}^3_+$ satisfying either

\[ \Phib_1 \big( \Xb^{*} \big) = \mathbf{0} \quad \text{and} \quad \Xb^{*} \in \mathbb{M}_s \]
or

\[ \Phib_2 \big( \Xb^{*} \big) = \mathbf{0} \quad \text{and} \quad \Xb^{*} \in \mathbb{M}_a \]
is called a \textit{virtual equilibrium} of the PWS system \eqref{sys}, \eqref{sys-pws}.

Let us also introduce for future use the following positive quantities:

\begin{equation}
\label{offspring-MF}
\N_M := \frac{\gamma r  \rho  \nu }{\mu  (\delta +\eta )} = \frac{r  \rho}{\mu} \cdot \frac{\gamma\nu }{\delta +\eta}, \qquad \N_{F}:= \dfrac{  (1-r) \rho \nu }{\delta  (\delta +\eta +\nu )} = \dfrac{  (1-r) \rho}{\delta} \cdot \dfrac{ \nu }{\delta +\eta +\nu}.
\end{equation}
These positive constants represent the basic offspring numbers related to the male and female psyllids. It is worthwhile to recall that the basic offspring number expresses a mean number of descendants produced by one individual during his/her lifespan. For males, $\N_M$ depends not only on the usual ratio $r \rho/\mu$ expressing an average number of eggs that later become males but also on the mating efficiency $\gamma$ of males and the relative availability for mating $\nu/(\delta + \nu)$ of the female psyllids. Similarly, for females, $\N_F$ depends not only on the usual ratio $(1-r) \rho/\delta$ expressing an average number of eggs that later become females but also on the mating frequency of female psyllids $\nu/(\delta +\eta +\nu)$.  Notably, the parameters $\nu, \eta$ related to the interchange between the compartments $A$ and $U$ are explicitly included in $\N_M$ and $\N_F$ meaning that the overall population size of adult insects strongly depends on the females' readiness for mating.

Now we proceed to identify the possible equilibria of smooth ODE systems \eqref{sys-F1} and \eqref{sys-F2}.

\subsection{Case 1: abundance of male psyllids}
\label{subsec-c1}

When $\gamma M > A $, the PWS system \eqref{sys} becomes \eqref{sys-F1} with $\Phib_1$ given by \eqref{F1}, and its equilibria are nonnegative solutions of the following algebraic system

\begin{subequations}
\label{fpsys-1}
\begin{align}[left = \empheqlbrace\,]
\label{fpsys-1-M}
0& = r \rho U e^{-\sigma (M+A+U)} - \mu M, \\[2mm]
\label{fpsys-1-A}
0& = (1-r) \rho U e^{-\sigma (M+A+U)} - \nu A +  \eta U - \delta A, \\[2mm]
\label{fpsys-1-U}
0& = \nu A - \eta U - \delta  U.
\end{align}
\end{subequations}
It is immediate to deduce that $\Eb_0=(0,0,0)$ is solution of \eqref{fpsys-1}. Then, we solve this system with $M$ and $A$ as unknowns and obtain

\begin{equation}
\label{fpsys1-MA}
M = \dfrac{r \delta   (\delta +\eta +\nu )}{(1-r) \mu  \nu } U, \qquad A = \dfrac{(\delta +\eta )}{\nu} U.
\end{equation}
Replacing these solutions in \eqref{fpsys-1-M}, we obtain

\[ r \rho U e^ {-\sigma (M+A+U)}=\frac{\delta  r (\delta +\eta +\nu )}{\nu  (1-r)}U \quad \Rightarrow \quad e^ {-\sigma (M+A+U)} = \frac{\delta  (\delta +\eta +\nu )}{ (1-r) \rho \nu  } = \frac{1}{\N_F} \]
according to the second relationship in \eqref{offspring-MF}. Thus, we deduce

\begin{equation}
\label{fpsys-sum1}
M+A+U= \frac{1}{\sigma} \ln \N_F>0
\end{equation}
meaning that a positive solution $\Eb^*_1:=\big( M_1^*, A_1^*, U_1^* \big)$ of \eqref{fpsys-1} exists if and only if $\N_F > 1.$ Further, by plugging the relationships \eqref{fpsys1-MA} into \eqref{fpsys-sum1} we obtain

\begin{equation}
\label{fp-sys1-aux}
 \frac{1}{\sigma} \ln \N_F = \left(\dfrac{r \delta   (\delta +\eta +\nu )}{(1-r) \mu  \nu } + \dfrac{(\delta +\eta )}{\nu} + 1 \right) U = \frac{(\delta +\eta +\nu ) \vartheta }{(1-r) \mu  \nu } U,
 \end{equation}
where

\begin{equation}
\label{vartheta}
\vartheta := (1-r) \mu +  r \delta
\end{equation}
denotes the so-called \textit{standardized mortality} of adult psyllids that, in effect, is the weighted mean mortality of both sex groups with the weights defined by their opposite-sex ratios. Finally, solving the equation \eqref{fp-sys1-aux} for $U$ and using it in \eqref{fpsys1-MA} we arrive to a strictly positive solution $\Eb^*_1:=\big( M_1^*, A_1^*, U_1^* \big)$ of \eqref{fpsys-1}:

\begin{subequations}
\label{E-1}
\begin{align}[left = \empheqlbrace\,]
\label{E-1M}
M_1^* &=r \frac{\delta}{\vartheta} \dfrac{1}{\sigma} \ln \N_F , \\	
\label{E-1A}
A_1^* &=  (1-r) \frac{\mu}{\vartheta} \frac{ \delta +\eta }{(\delta +\eta +\nu)} \dfrac{1}{\sigma} \ln \N_F , \\
\label{E-1U}
U_1^* &=  (1-r) \frac{\mu}{\vartheta}  \frac{\nu}{( \delta +\eta +\nu)} \dfrac{1}{\sigma} \ln \N_F .
\end{align}
\end{subequations}
Note also that the total number of insects at equilibrium $\Eb_1^*$ verifies \eqref{fpsys-sum1}, and its coordinates explicitly include the standardized mortality ratios ($r\delta/\vartheta$ and $(1-r)\mu/\vartheta$) related to opposite sex.

Thus, we conclude that the smooth system \eqref{sys-F1} has two possible equilibria: the trivial equilibrium $\Eb_0=(0,0,0)$ that exists for any positive value of $\N_F$ (defined by \eqref{offspring-MF}), and the strictly positive one $\Eb^*_1=\big( M_1^*, A_1^*, U_1^* \big)$ defined by \eqref{E-1} that exists if and only if $\N_F >1.$ The following result establishes the stability properties of $\Eb_0$ and $\Eb_1^*$.

\begin{proposition}
\label{prop2}
\begin{itemize}
\item
Assume $\N_F<1$. Then $\Eb_0$ is locally asymptotically stable (LAS).
\item
Assume $\N_F>1$. Then $\Eb_1^*$ is LAS and $\Eb_0$ is unstable; however, there always exists a trajectory converging to $\Eb_0$ meaning that $\Eb_0$ is not a repeller.
\end{itemize}
\end{proposition}
\begin{proof}
See Appendix A, page \pageref{appendixA}.
\end{proof}

\subsection{Case 2: scarcity of male psyllids}
\label{subsec-c2}

When $\gamma M < A $, the PWS system \eqref{sys} becomes \eqref{sys-F2} with $\Phib_2$ given by \eqref{F2}, and its equilibria are nonnegative solutions of the following algebraic system

\begin{subequations}
\label{fpsys-2}
\begin{align}[left = \empheqlbrace\,]
\label{fpsys-2-M}
0& = r \rho U e^{-\sigma (M+A+U)} - \mu M, \\[2mm]
\label{fpsys-2-A}
0& = (1-r) \rho U e^{-\sigma (M+A+U)} - \gamma \nu M  +  \eta U - \delta A, \\[2mm]
\label{fpsys-2-U}
0& = \gamma \nu M - \eta U - \delta  U.
\end{align}
\end{subequations}
It is immediate to deduce that $\Eb_0=(0,0,0)$ is solution of \eqref{fpsys-2}. We also set

\begin{equation}
\label{theta_M}
\theta_M := \frac{  (1-r) \mu (\delta +\eta )}{\gamma r \delta  \nu}.
\end{equation}

Then, we solve the nonlinear system \eqref{fpsys-2} with $M$ and $A$ as unknowns and obtain

\begin{equation}
\label{fpsys2-MA}
M  = \dfrac{ (\delta +\eta )}{\gamma  \nu } U, \qquad A =  \left(\frac{\mu  (1-r) (\delta +\eta )}{\gamma  \delta  \nu  r}-1\right) U= \big(\theta_M-1 \big) U.
\end{equation}
Thus, $A>0$ whenever $\theta_M > 1.$ It is interesting to notice that

\[ \theta_M=\dfrac{1}{\N_M} \dfrac{(1-r)\rho}{\delta}, \]
such that having $\N_M>1$ (cf. the first relationship in \eqref{offspring-MF}), we need

\begin{equation}
\label{thetaM2}
    \dfrac{(1-r)\rho}{\delta}> \N_M,
\end{equation}
in order to assure that $\theta_M>1$. In fact, \eqref{thetaM2} means that the average number of eggs that further become females has to be larger than the mean number of male descendants produced by one male individual all along his lifespan. In other words, condition $\theta_M>1$ can be replaced by \eqref{thetaM2}.

Direct substitution of \eqref{fpsys2-MA} into \eqref{fpsys-2-M} renders

\[ r \rho U e^{-\sigma (A+M+U)} =\dfrac{\mu  (\delta +\eta)}{\gamma \nu } U \quad \Rightarrow \quad  e^{-\sigma (A+M+U)} =\dfrac{\mu  (\delta +\eta)}{\gamma r \rho  \nu } = \frac{1}{\N_M}  \]
leading to

\begin{equation}
\label{fpsys-sum2}
 A+M+U=\dfrac{1}{\sigma}\ln{\N_M} > 0,
\end{equation}
and meaning that a positive solution $\Eb_2^*= \big( M_2^*, A_2^*, U_2^* \big)$ of \eqref{fpsys-2} exists if and only if $\N_M > 1$ and $\theta_M > 1$, where $\N_M$ and $\theta_M$ are given by \eqref{offspring-MF} and \eqref{theta_M}, respectively. Further, by plugging the relationships \eqref{fpsys2-MA} into \eqref{fpsys-sum2} we obtain

\begin{equation}
\label{fp-sys2-aux}
\frac{1}{\sigma}\ln \N_M = M+A+U=\left( \frac{\eta+\delta}{\nu\gamma} + \theta_M - 1 + 1 \right) U= \frac{\gamma \nu \theta_M + \eta + \delta}{\gamma \nu} U.
\end{equation}
Finally, solving equation \eqref{fp-sys2-aux} for $U$ and using it in \eqref{fpsys2-MA} we arrive to a strictly positive solution $\Eb_2^*:= \big( M_2^*, A_2^*, U_2^* \big)$ of \eqref{fpsys-2}:

\begin{subequations}
\label{equilib2}
\begin{align}[left = \empheqlbrace\,]
\label{equilib2-M}
M_2^* &= \frac{\delta +\eta}{\gamma  \nu \theta_M +  \eta+\delta}\frac{1}{\sigma}\ln{\N_M}, \\
\label{equilib2-A}
A_2^* &= \frac{\gamma \nu \big( \theta_M-1 \big)}{\gamma  \nu \theta_M +  \eta+\delta}\frac{1}{\sigma}\ln{\N_M}\\
\label{equilib2-U}
U_2^* &= \frac{\gamma \nu }{\gamma  \nu \theta_M +  \eta+\delta} \frac{1}{\sigma}\ln{\N_M},
\end{align}
\end{subequations}

\begin{remark}
One can also obtain the expressions for coordinates of $\Eb_2^*$ in terms of the standardized mortality $\vartheta$ defined by \eqref{vartheta}. Notably, the denominator of all expressions included in \eqref{equilib2} verifies

\[ \gamma  \nu \theta_M +  \eta+\delta = \gamma \nu \frac{(1-r) \mu (\delta + \eta)}{\gamma  r \delta  \nu} + (\delta + \eta) = (\delta + \eta) \left( \frac{(1-r) \mu }{ r \delta} + 1 \right) =  (\delta + \eta) \frac{\vartheta}{r \delta}. \]
Furthermore,

\[ \gamma \nu \big( \theta_M - 1 \big) = \frac{(1-r) \mu (\delta + \eta)}{r \delta} - \gamma \nu = \frac{\mu}{r} (\delta + \eta) \left( \frac{1-r}{\delta} - \frac{\gamma r \nu}{\mu (\delta + \eta)} \right)= \frac{\mu}{r \rho} (\delta + \eta) \left( \frac{(1-r) \rho}{\delta} - \N_M \right).
\]
Using the above relationships in combination with \eqref{equilib2}, we can obtain an alternative form of \eqref{equilib2}:

\begin{subequations}
\label{E-2}
\begin{align}[left = \empheqlbrace\,]
\label{E-2M}
M_2^* & = r \frac{\delta}{\vartheta} \frac{1}{\sigma} \ln \N_M, \\
\label{E-2A}
A_2^* & = \frac{\mu}{\vartheta} \frac{\delta}{\rho} \left( \frac{(1-r) \rho}{\delta} - \N_M \right) \frac{1}{\sigma} \ln \N_M, \\
\label{E-2U}
U_2^* & =\frac{\mu}{\vartheta} \frac{\delta}{\rho} \N_M \frac{1}{\sigma}\ln{\N_M},
\end{align}
\end{subequations}
This alternative form of $\Eb_2^*$ makes visible the necessity of the condition \eqref{thetaM2} for existence of $\Eb_2^*$ along with $\N_M >1.$
\end{remark}

Note also that for both forms of $\Eb_2^*$ (\eqref{equilib2} and \eqref{E-2}),  the total number of insects at equilibrium $\Eb_2^*$ verifies \eqref{fpsys-sum2}.

Thus, we conclude that the smooth system \eqref{sys-F2} has two possible equilibria: the trivial equilibrium $\Eb_0=(0,0,0)$ that exists for any positive value of $\N_M$ (defined by \eqref{offspring-MF}), and the strictly positive one $\Eb^*_2=\big( M_2^*, A_2^*, U_2^* \big)$ defined by \eqref{equilib2} or \eqref{E-2} that exists if and only if $\N_M >1$ and $\theta_M >1 $, that is, if the condition \eqref{thetaM2} holds. The following result establishes the stability properties of $\Eb_0$ and $\Eb_2^*$.

\begin{proposition}
\label{prop3}
\begin{itemize}
\item
Assume $\N_M<1$. Then $\Eb_0$ is locally asymptotically stable (LAS).
\item
Assume $\N_M>1$ and $\theta_M > 1$ meaning that \eqref{thetaM2} holds. Then $\Eb_2^*$ is LAS and $\Eb_0$ is unstable; however, there always exists a trajectory converging to $\Eb_0$ meaning that $\Eb_0$ is not a repeller.
\end{itemize}
\end{proposition}
\begin{proof}
See Appendix A, page \pageref{appendixA}.
\end{proof}

\subsection{Stability appraisal for the PWS system \eqref{sys}}
\label{subsec-pws}

First, we note that $\Eb_0 \in \mathcal{P}_s.$ Let us now determine the position of $\Eb_1^*$ and $\Eb_2^*$ in $\mathbb{R}^3_+$ with respect to the switching plane $\mathcal{P}_s$. Clearly,   $\Eb_1^* \in \mathcal{P}_s$ if and only if $\gamma M_1^* = A_1^*$, that is,

\[ \gamma r \frac{\delta}{\vartheta} \dfrac{1}{\sigma} \ln{\N_F} = (1-r) \frac{\mu}{\vartheta} \frac{ \delta +\eta }{(\delta +\eta +\nu)} \dfrac{1}{\sigma} \ln{\N_F} \qquad \text{or} \qquad  \gamma r \delta = (1-r) \mu \frac{ \delta +\eta }{(\delta +\eta +\nu)}. \]

By multiplying both sides of the last relationship by $\dfrac{\rho \nu}{\mu \delta (\delta + \eta)} >0$ we arrive to

\[ \frac{\gamma r  \rho  \nu }{\mu  (\delta +\eta )} = \dfrac{(1-r) \rho \nu }{\delta  (\delta +\eta +\nu )} \qquad \Leftrightarrow \qquad \N_M = \N_F. \]
Thus, $\Eb_1^* \in \mathcal{P}_s$ if and only if $\N_M = \N_F > 1$. Furthermore, it is easy to deduce that

\[ \Eb_1^* \in \mathbb{M}_a \quad \Leftrightarrow \quad \N_M > \N_F > 1 \quad \text{and} \quad \Eb_1^* \in \mathbb{M}_s \quad \Leftrightarrow \quad \N_M < \N_F, \; \N_F > 1. \]
The above expressions imply that $\Eb_1^*$ is a \textit{regular equilibrium} of the original PWS system \eqref{sys} when $\N_M > \N_F > 1$, and $\Eb_1^*$ is a \textit{virtual equilibrium} of \eqref{sys} when $\N_M < \N_F$ and $\N_F > 1$.

\begin{figure}[t]
\centering
	\includegraphics[width=10cm]{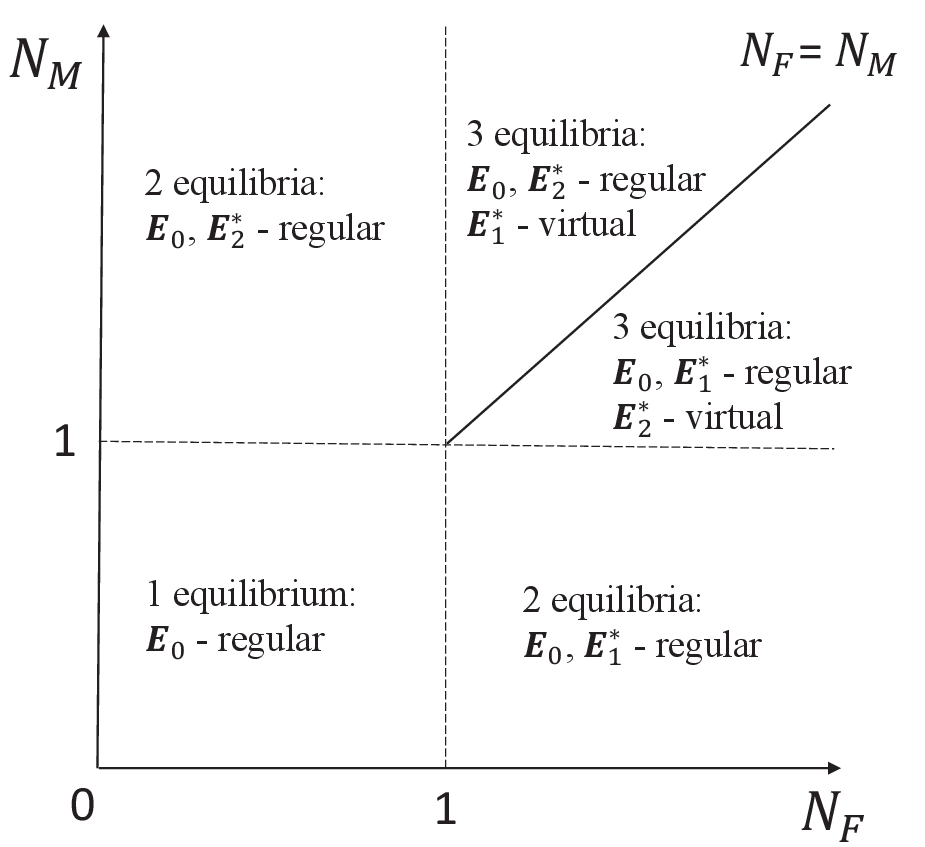}
\caption{ Regular and virtual equilibria of the PWS system \eqref{sys-pws}-\eqref{F2} according to the values of $\N_F$ and $\N_M$ \label{fig-equi}}
\end{figure}

On the other hand, $\Eb_2^* \in \mathcal{P}_s$ if and only if $\gamma M_2^* = A_2^*$, that is,

\[ \frac{\gamma (\delta +\eta)}{\gamma \nu \theta_M + \eta+\delta}\frac{1}{\sigma}\ln{\N_M} = \frac{\nu \gamma \big( \theta_M-1 \big)}{\gamma \nu \theta_M + \eta+\delta} \frac{1}{\sigma} \ln \N_M \qquad \text{or} \qquad \delta + \eta + \nu = \frac{\nu(1-r) \mu (\delta +\eta )}{\gamma r \delta  \nu}. \]
By multiplying both sides of the last relationship by $\dfrac{\gamma r \rho \nu}{\mu (\delta + \eta) (\delta + \eta + \nu)} >0$ we arrive to

\[ \frac{\gamma r  \rho  \nu }{\mu  (\delta +\eta )} = \dfrac{(1-r) \rho \nu }{\delta (\delta +\eta +\nu )} \qquad \Leftrightarrow \qquad \N_M = \N_F. \]
Thus, $\Eb_2^* \in \mathcal{P}_s$ if and only if $\N_M = \N_F > 1$. Furthermore, it is easy to deduce that

\[ \Eb_2^* \in \mathbb{M}_s \quad \Leftrightarrow \quad \N_F > \N_M > 1 \quad \text{and} \quad \Eb_2^* \in \mathbb{M}_a \quad \Leftrightarrow \quad \N_M > \N_F, \; \N_M > 1. \]
The above expressions imply that $\Eb_2^*$ is a \textit{regular equilibrium} of the original PWS system \eqref{sys} when $\N_F > \N_M > 1$, and $\Eb_2^*$ is a \textit{virtual equilibrium} of \eqref{sys} when $\N_M > \N_F$ and $\N_M > 1$. Figure \ref{fig-equi}  schematically displays the regular and virtual equilibria the PWS system \eqref{sys-pws}-\eqref{F2} may possess according to the values of the basic offspring numbers $\N_F$ and $\N_M$.

From the foregoing rationale, we can also conclude that $\N_M = \N_F > 1$ implies that $\Eb_1^* = \Eb_2^* \in \mathcal{P}_s$ meaning that both positive equilibria collide and coalesce on the switching plane $\mathcal{P}_s$. The latter can be checked by comparing the components of $\Eb_1^*$ and $\Eb_2^*$ using their forms given by \eqref{E-1} and \eqref{E-2} when $\N_M = \N_F > 1$.

\section{Pest control by pheromone traps}
\label{sec-discussion}

Female psyllids available for mating (class $ A $) emit sex pheromones that attract male insects over a long distance \cite{Li2020}.

Sex pheromone traps offer an alternative to traditional pesticides and can be considered an eco-friendly component of integrated pest control. First,   pheromone traps can be used for monitoring pest insects to determine whether additional control measures are needed. Second, pheromone traps can be set up as a lure to perform control of pest populations. In this case, sticky pheromone traps emitting large quantities of sex pheromones may serve one of the following two purposes or both of them:
\begin{enumerate}
\item
Attraction and mass trapping of male insects, followed by their direct removal (male killing).
\item
Mating disruption for decreasing the fecundity of females (offspring reduction).
\end{enumerate}

The ACP population dynamics model \eqref{sys} proposed in Section \ref{sec-model} can be adapted to include the two control actions mentioned above. Let $A_p >0$ express the ``strength of lure''. Knowing the average amount of sex pheromones emitted by one female psyllid \cite{Mann2013}, the external parameter $A_p$ can be expressed in terms of the number of ``false'' female psyllids available for mating. Then, the total number of male-seeking females (both natural and false) is expressed by $(A_p + A)$ \cite{Anguelov2017,Barclay1983}. Furthermore, a female-seeking male has the probability $\dfrac{A}{A_p + A}$ of being attracted to a wild (natural) female and the probability $\dfrac{A_p}{A_p + A}$ of being attracted to the pheromone traps. Let $\alpha \in [0,1]$ denote the capture or killing rate of males attracted to a pheromone trap. Then, by setting $\alpha=1$, it is modeled that all males approaching or entering the trap are killed, while $\alpha =0$ models that none of them will be killed when approaching or entering the trap. Notably, by setting $\alpha=0$ and $A_p=0$, the original model \eqref{sys} can be immediately recovered.

Using two additional parameters ($A_p$ and $\alpha$) defined above, we can now formulate the modified version of the model \eqref{sys} that accounts for mating disruption and male-killing effect induced by the pheromone traps:

\begin{subequations}
\label{syscon}
\begin{align}[left = \empheqlbrace\,]
\label{syscon-M}
\frac{d M}{dt}& = r \rho U e^{-\sigma (M+A+U)} - \alpha \frac{A_p}{A_p + A} M - \mu M \\
\label{syscon-A}
\frac{d A}{dt}& = (1-r) \rho U e^{-\sigma (M+A+U)} - \nu \min \left\{ \frac{\gamma M}{A_p + A}, 1  \right\} A +  \eta U - \delta A \\
\label{syscon-U}
\frac{d U}{dt}& = \nu \min \left\{ \frac{\gamma M}{A_p + A}, 1  \right\} A - \eta U - \delta  U
\end{align}
\end{subequations}

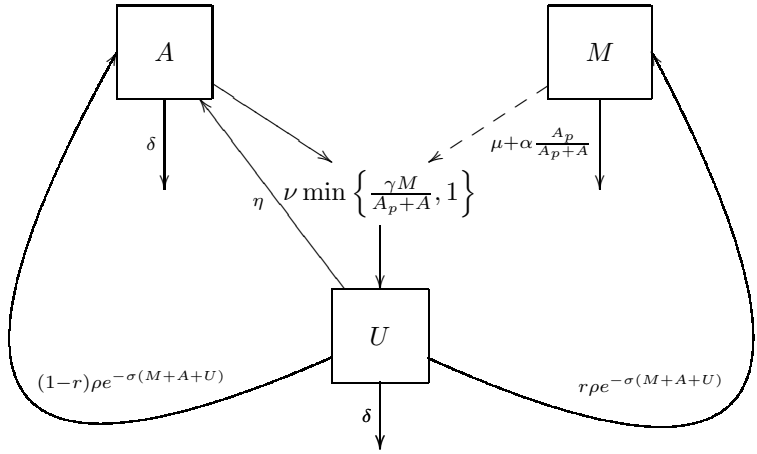
\begin{figure}[t]
\centering
	\pagestyle{empty}
	\thispagestyle{empty}
	\hspace{-2.5cm}
  \centering
	\xymatrix{
	*+<1cm>[F]{A}\ar@{->}[dr] \ar[d]_-{\delta} & & *+<1cm>[F]{M}\ar@{-->}[dl]  \ar[d]_-{\mu +\alpha \frac{A_p}{A_p + A}}\\
	&{\nu \min \left\{ \frac{\gamma M}{A_p + A}, 1  \right\}}\ar[d] & \\
	& *+<1cm>[F]{U}\ar@{->}[uul]^{\eta}  \ar[d]_-{\delta} \ar@[mosco]@/^4cm/[luu]+L_{(1-r) \rho e^{-\sigma (M+A+U)}} \ar[d]_-{\delta} \ar@[mosco]@/^-4cm/[uur]+R^{r \rho e^{-\sigma (M+A+U)}} & \\
	& & \\
	}
\caption{Flow diagram of the population dynamics of \textit{D. citri} with pheromone traps \eqref{syscon} \label{fig-2}}
\end{figure}

The initial conditions for this model are the same as \eqref{icon}, and Figure \ref{fig-2} provides the flow diagram of the model \eqref{syscon}.

Similarly to the original dynamical system \eqref{sys}, the model \eqref{syscon} is a PWS system that features two external parameters $\alpha \in [0,1]$ and $A_p \geq 0$, one of which directly affects its switching plane

\[ \widetilde{\mathcal{P}}_s (A_p) := \Big\{ (M,A,U) \in \mathbb{R}^3_{+} : \ \gamma M = A + A_p \Big\}. \]

\begin{figure}[t]
\centering
	\includegraphics[width=13cm]{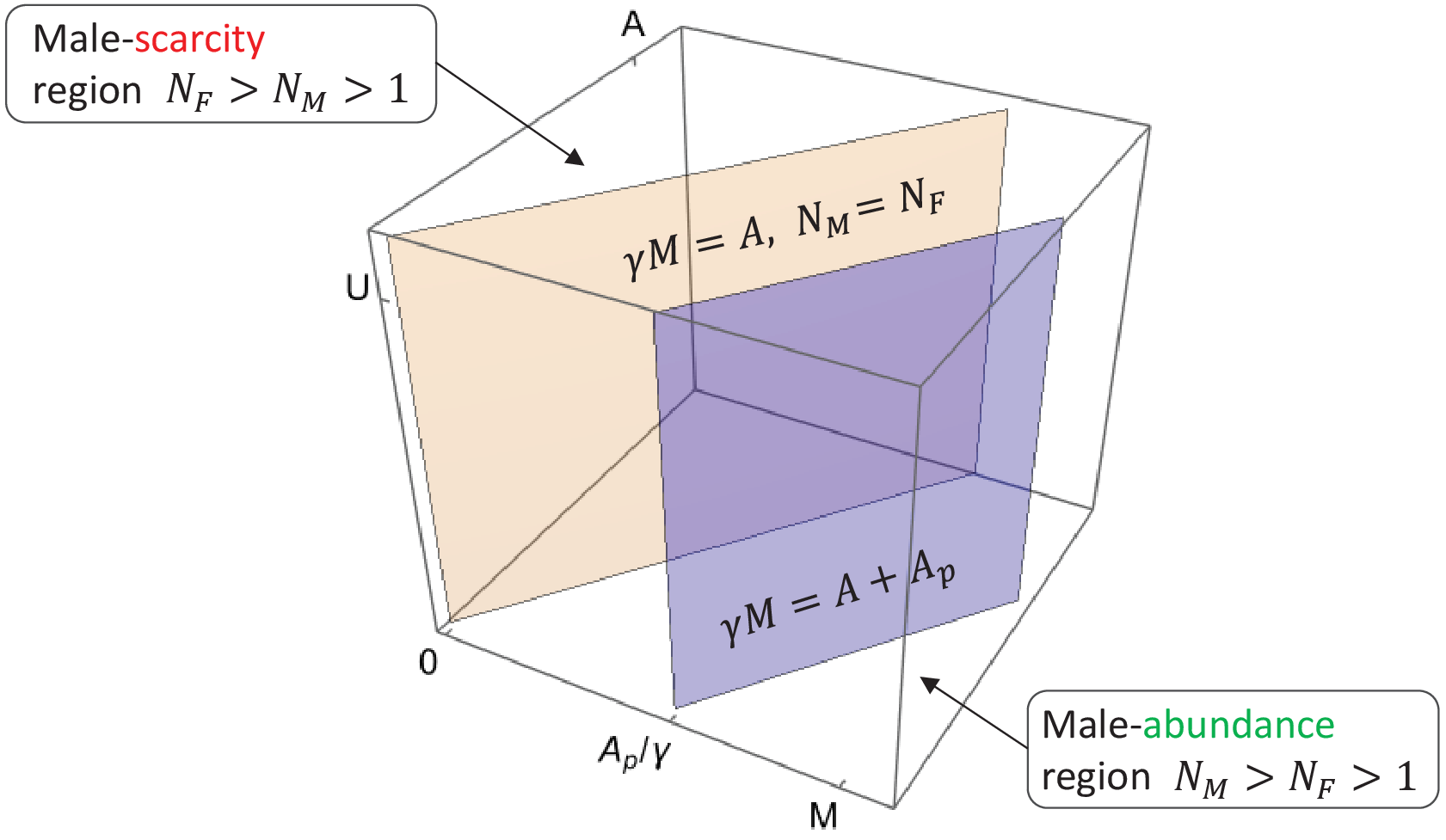}
\caption{ Changes in the male-scarcity and male-abundance regions induced by $A_p > 0$ \label{fig-planes}}
\end{figure}

The geometric role of the parameter $A_p >0$ is illustrated in Figure \ref{fig-planes}: for larger values of $A_p>0$, the male-scarcity region becomes more extensive, and the switching plane $ \widetilde{\mathcal{P}}_s (A_p)$ moves farther away from the origin.

Existence and uniqueness of a piecewise smooth solution to the system \eqref{syscon} for any nonnegative initial conditions \eqref{icon}, as well as the nonnegativity and boundedness of the system's trajectories in $\mathbb{R}_+^3$ can be proved using the arguments presented in the proof of Proposition \ref{prop1}.

To analyze the behavior of solutions of the PWS system \eqref{syscon}, we split it, following the approach in \cite{Anguelov2017}, into two smooth systems, each featuring two external parameters. For each smooth system, we propose two operational control modes referred to as open-loop and closed-loop control approaches in the sequel. It is worthwhile to recall here that open-loop control operates based on a predefined sequence of actions. In contrast, the closed-loop control approach is more adaptive since its actions can respond to changes in the system behavior.

\subsection{Open-loop control approach}
\label{subsec-open}

In the male abundance region $\mathbb{M}_a$, where it holds $\gamma M>A+A_p$, the dynamical system \eqref{syscon} takes the following form:

\begin{subequations}
\label{syscon_abundance_open_loop}
\begin{align}[left = \empheqlbrace\,]
\label{syscon-M1}
\frac{d M}{dt}& = r \rho U e^{-\sigma (M+A+U)} - \alpha \frac{A_p}{A_p + A} M - \mu M \\
\label{syscon-A1}
\frac{d A}{dt}& = (1-r) \rho U e^{-\sigma (M+A+U)} - \nu A +  \eta U - \delta A \\
\label{syscon-U1}
\frac{d U}{dt}& = \nu  A - \eta U - \delta  U
\end{align}
\end{subequations}
and in the male-scarcity region $\mathbb{M}_s$, where it holds $\gamma M<A+A_p$, the dynamical system \eqref{syscon} becomes:

\begin{subequations}
\label{syscon_scarse_open_loop}
\begin{align}[left = \empheqlbrace\,]
\label{syscon-M2}
\frac{d M}{dt}& = r \rho U e^{-\sigma (M+A+U)} - \alpha \frac{A_p}{A_p + A} M - \mu M \\
\label{syscon-A2}
\frac{d A}{dt}& = (1-r) \rho U e^{-\sigma (M+A+U)} - \nu \frac{\gamma M}{A_p + A} A +  \eta U - \delta A \\
\label{syscon-U2}
\frac{d U}{dt}& = \nu \frac{\gamma M}{A_p + A} A- \eta U - \delta  U
\end{align}
\end{subequations}
Following the same rationale as in Subsection \ref{subsec-c1}, it is straightforward to show that system \eqref{syscon_abundance_open_loop} admits two equilibria, $\Eb_0$, and $\Eb_{1,P}^*=(M^*_{1,P},A^*_{1,P},U^*_{1,P})$, where $A^*_{1,P}$ is the positive root of the quadratic equation

\[
\left( \! \dfrac{\delta + \nu + \eta}{\delta + \eta} \! \right) \! \left( \! \dfrac{(1-r) \mu + \delta r}{1-r} \! \right) A^{2} + \left[ \! \left( \! \dfrac{\delta r}{1-r} + \left( \mu + \alpha \right) \! \right) \! \left( \! \dfrac{\delta + \nu + \eta}{\delta + \eta} \! \right) A_{p} - \dfrac{\mu}{\sigma} \ln \N_{F} \right] A - \left(\mu + \alpha \right) A_{p} \dfrac{1}{\sigma} \ln \N_{F}=0,
\]
or

\[
\left( \dfrac{\delta + \nu + \eta}{\delta + \eta}\right) \left( \dfrac{\vartheta}{1-r} \right) A^{2} + \left[ \left( \dfrac{\vartheta}{1-r} + \alpha \right) \left( \dfrac{\delta + \nu + \eta}{\delta + \eta}\right) A_{p} - \dfrac{\mu}{\sigma} \ln \N_{F} \right] A - \left( \mu+\alpha \right) A_{p} \dfrac{1}{\sigma} \ln \N_{F}=0.
\]
Notably, the above equation has only one positive real root because its discriminant is positive, the branches of the corresponding parabola are directed upwards,  while its cross with the vertical axis is negative. By setting $\Delta$ as the discriminant  of the above quadratic equation, we deduce

\begin{align*}
  A^*_{1,P} &= \dfrac{1}{2 \left( \dfrac{\delta + \nu + \eta}{\delta + \eta} \right) \left( \dfrac{(1-r)\mu + \delta r}{1-r} \right)} \left[ \dfrac{\mu}{\sigma} \ln \N_{F} - \left( \dfrac{r}{1-r} + \left( \mu + \alpha \right) \right) \left( \dfrac{\delta + \nu + \eta}{\delta + \eta} \right) A_{p} + \sqrt{\Delta} \right],   \\[2mm]
  U^*_{1,P} &=\dfrac{\nu}{\delta + \eta} A^*_{1,P}, \\[2mm]
  M^*_{1,P} &= \dfrac{r \rho}{\mu + \alpha \dfrac{A_{p}}{ A^*_{1,P} + A_{p}}} \: \dfrac{1}{\N_{F}} \: \dfrac{\nu}{\delta + \eta} A^*_{1,P}.
\end{align*}
Note also that $M^*_{1,P} + A^*_{1,P} + U^*_{1,P}=\dfrac{1}{\sigma} \ln \N_F$  meaning that, in the male abundance case, there is no impact on the total population. Of course, when $A_p=0$, we recover the positive equilibrium $\Eb_1^*$ given by \eqref{E-1} in Subsection \ref{subsec-c1}.

Last but not least, following the same methodology as in Appendix A, it is straightforward to show that for system \eqref{syscon_abundance_open_loop}, we have the following result.
\begin{proposition}
Assume $A_p\geq 0$.
\label{prop2_open_loop}
\begin{itemize}
\item
If $\N_F<1$, then $\Eb_0$ is locally asymptotically stable (LAS).
\item
If $\N_F>1$, then $\Eb_{1,P}^*$ is LAS and $\Eb_0$ is unstable.
\end{itemize}
\end{proposition}

Let us now focus on system \eqref{syscon_scarse_open_loop}. Like in Subsection \ref{subsec-c2}, it is straightforward to show that $\Eb_0$ is still an equilibrium. However, mating disrupting entails the Allee effect when the system parameters correspond to the male-scarcity case.  Namely, $\Eb_0$ may remain locally asymptotically stable as long as \textit{a sufficient amount of} pheromones are released and regardless of $\N_F$. The latter can be shown via direct computation of the Jacobian matrix at equilibrium $\Eb_0$. This property is essential because it allows the use of a small amount of pheromones to control a non-established population and, from the field application point of view, to derive a massive and small releases strategy, like in the Sterile Insect Technique \cite{Anguelov2020}. Notwithstanding, showing the existence of at least one positive equilibrium for system \eqref{syscon_scarse_open_loop} is a bit more complicated than in Subsection \ref{subsec-c2}.

\begin{proposition}
\label{existence_equilibrium_scarse_pheromones}
    Assume that $A_p\geq 0$ and $\theta_M>1$, where $\theta_M$ is defined by \eqref{theta_M}. There exists a threshold quantity $A_p^{crit}>0$ such that
    \begin{itemize}
        \item When $0<A_p<A_p^{crit}$, system \eqref{syscon_scarse_open_loop} admits 2 positive equilibria, $\Eb_{1,p}$ and $\Eb_{2,p}$.
        \item When $A_p=A_p^{crit}$, system \eqref{syscon_scarse_open_loop} admits only one positive equilibrium $\Eb_{*,p}$.
        \item When $A_p>A_p^{crit}$, system \eqref{syscon_scarse_open_loop} has no positive equilibrium.
    \end{itemize}
\end{proposition}
\begin{proof}
    See Appendix B, page \pageref{appendixB}.
\end{proof}

From Proposition \ref{existence_equilibrium_scarse_pheromones}, when $A_p>A_p^{crit}$, we deduce that the only equilibrium is $\bf{E}_0$, which is always LAS.
\begin{remark}
    The previous result shows that massive releases of pheromones can be used to suppress or eradicate the ACP population. However, the emission of pheromones in large quantities is not always necessary as long as we can estimate the wild population during the intervention. That is why the closed-loop control approach can help derive some strategies relying on a minimal amount of pheromones. Subsection \ref{subsec-closed} addresses this issue in more detail.
\end{remark}

It is also possible to show that $\bf{E}_0$ is not only LAS but also GAS (\textit{globally asymptotically stable}) when $A_p$ is sufficiently large. Following \cite{Anguelov2017}, it is easy to check that the right-hand side of system \eqref{syscon_scarse_open_loop} is not quasi-monotone because of the term $- \nu \dfrac{\gamma M}{A_p + A} A$. However, it is possible to consider an auxiliary system that is monotone cooperative and provides an upper-bound for the solution of system \eqref{syscon_scarse_open_loop} by removing $- \nu \dfrac{\gamma M}{A_p + A} A$ and also some exponential terms. The auxiliary system becomes

\begin{subequations}
\label{auxiliary_scarse_open_loop}
\begin{align}[left = \empheqlbrace\,]
\label{syscon-Ma}
\frac{d M}{dt}& = r \rho U e^{-\sigma M} - \alpha \frac{A_p}{A_p + A} M - \mu M \\
\label{syscon-Aa}
\frac{d A}{dt}& = (1-r) \rho U  +  \eta U - \delta A \\
\label{syscon-Ua}
\frac{d U}{dt}& = \nu \frac{\gamma M}{A_p + A} A- \eta U - \delta  U.
\end{align}
\end{subequations}

We can derive the following result.
\begin{theorem}
\label{theo_GAS_scarse_aux}
\begin{itemize}
    \item[(a)] System \eqref{auxiliary_scarse_open_loop} defines a positive dynamical system on $\mathbb{R}^3_+$.
    \item[(b)] Equilibrium $\Eb_0=\bf{0}$ of the system \eqref{auxiliary_scarse_open_loop} is always LAS.
    \item[(c)] There exists $\tilde{A}_p^{crit}>0$ such that
    \begin{itemize}
        \item If $A_p>\tilde{A}_p^{crit}$, then $\Eb_0=\bf{0}$ is GAS on $\mathbb{R}^3_+$.
        \item If $0<A_p<\tilde{A}_p^{crit}$, system \eqref{auxiliary_scarse_open_loop} admits two positive equilibria $\tilde{\Eb}_{1,p}$ and $\tilde{\Eb}_{2,p}$ such that $\tilde{\Eb}_{1,p}<\tilde{\Eb}_{2,p}$. Moreover, the basin of attraction of $\Eb_0$ contains the set $\Big\{ \Xb \in \mathbb{R}_+^3: \ \Eb_0 \leq  \Xb < \tilde{\Eb}_{1,p} \Big\}$, and the basin of attraction of $\tilde{\Eb}_{2,p}$ contains the set $\Big\{ \Xb \in \mathbb{R}_+^3: \ \Xb \geq  \tilde{\Eb}_{2,p} \Big\}$.
    \end{itemize}
\end{itemize}
\end{theorem}
\begin{proof}
    See Appendix C, page \pageref{appendixC}
\end{proof}
Finally, by comparison, any solution of \eqref{auxiliary_scarse_open_loop} is an upper bound for the solution of \eqref{syscon_scarse_open_loop} with the same initial point. This implies that the basin of attraction of $\Eb_0$ as an equilibrium of \eqref{syscon_scarse_open_loop} contains the sets given in Theorem \ref{theo_GAS_scarse_aux} (c). Thus we deduce the following statement.
\begin{theorem}
\label{theo_GAS_scarse_model}
    Let $A_p>0$. Then, the following hold for model \eqref{syscon_scarse_open_loop}:
    \begin{itemize}
        \item[(a)]  If $0<A_p<\tilde{A}_p^{crit}$, the basin of attraction of $\Eb_0$ contains the set $\Big\{ \Xb \in \mathbb{R}_+^3: \ \Eb_0 \leq \Xb < \tilde{\Eb}_{1,p} \Big\}$.
        \item[(b)] If $A_p >\tilde{A}_p^{crit}$, then $\Eb_0=\bf{0}$ is GAS on $\mathbb{R}^3_+$.
    \end{itemize}
\end{theorem}

Theorem \ref{theo_GAS_scarse_model} provides the following helpful information from the practical perspective:
\begin{itemize}
    \item As long as $A_p<\tilde{A}_p^{crit}$, only an invading or non-established population can be controlled.
    \item An established population can only be eliminated when  $A_p>\tilde{A}_p^{crit}$.
\end{itemize}

\subsection{Closed-loop control approach}
\label{subsec-closed}

Let us now assume that the total amount of sex pheromones expressed in terms of ``false'' females is proportional to the number of natural female psyllids $A(t)$ seeking for mating that is,

\begin{equation}
\label{gain}
 A_p =  k \times A(t),
\end{equation}
where the constant $k > 0$ defines the ``gain'' of feedback. Then, in the male-abundance region $\gamma M > A + A_p =(k+1)A$, the dynamical system \eqref{syscon} takes the following form:

\begin{subequations}
\label{sys-abun}
\begin{align}[left = \empheqlbrace\,]
\dfrac{dM}{dt} & = r\rho U e^{-\sigma(M+A+U)} - \alpha \dfrac{k}{k+1} M -\mu M \\
\dfrac{dA}{dt} & =  (1-r)\rho U e^{-\sigma(M+A+U)} - \nu A + \eta U - \delta A \\
\dfrac{dU}{dt} & =  \nu A - \eta U - \delta U
\end{align}
\end{subequations}

Alternatively, in the male-scarcity region $\gamma M < A + A_p =(k+1)A$, the dynamical system \eqref{syscon} becomes

\begin{subequations}
\label{sys-scar}
\begin{align}[left = \empheqlbrace\,]
\dfrac{dM}{dt} & = r\rho U e^{-\sigma(M+A+U)} - \alpha \dfrac{k}{k+1} M -\mu M \\
\dfrac{dA}{dt} & =  (1-r)\rho U e^{-\sigma(M+A+U)} - \dfrac{\gamma \nu}{k+1} M + \eta U - \delta A \\
\dfrac{dU}{dt} & =  \dfrac{\gamma \nu}{k+1} M - \eta U - \delta U
\end{align}
\end{subequations}

Let us now derive, for both systems \eqref{sys-abun} and \eqref{sys-scar}, the conditions that define either the permanence or extinction of the ACP population under the feedback \eqref{gain} and in the presence of male-killing traps. To do so, we will employ the next-generation approach initially derived for epidemiological systems \cite{Driessche2002} and later adapted to more general population dynamics models, see for instance \cite{Barril2018,Barril2021}. This approach consists in determining the spectral radius of the next-generation matrix evaluated in the trivial equilibrium of the population dynamics model. To construct the next-generation matrix, the right-hand side of the dynamical system is written in the form

\[ \dfrac{d \Xb}{dt} = \Gb(\Xb):=  \mathcal{F}(\Xb) - \mathcal{V}(\Xb), \]
where the vector $\mathcal{F}(\Xb)$ gathers all terms dealing with the emergence of new individuals, while the vector $\mathcal{V}(\Xb)$ contains the transition and mortality terms. For the systems \eqref{sys-abun} and \eqref{sys-scar}, the vector $\mathcal{F}(\Xb)$ is same, that is,

\[ \mathcal{F}(\Xb) = \begin{pmatrix}
r \rho \: U e^{-\sigma (M + A + U)} \\[2mm]
(1-r) \rho \: U e^{-\sigma (M + A + U)} \\[2mm]
0
\end{pmatrix}, \qquad \Xb = \begin{pmatrix} M \\  A \\ U \end{pmatrix},
\]
while $\mathcal{V}(\Xb)$ takes different forms:

\[ \text{System \eqref{sys-abun} } \Rightarrow \mathcal{V}(\Xb) \!=\! \mathcal{V}_1 (\Xb) \!=\!
\begin{pmatrix} \dfrac{\alpha k}{k+1} M \!+\! \mu M \\[3mm]
(\nu \!+\! \delta) A \!-\! \eta U \\[3mm]
(\eta \!+\! \delta)U \!-\! \nu A
\end{pmatrix}\!, \quad \text{System \eqref{sys-scar} } \Rightarrow \mathcal{V}(\Xb) \!=\! \mathcal{V}_2 (\Xb) \!=\!
\begin{pmatrix} \dfrac{\alpha k}{k+1} M \!+\! \mu M \\[3mm]
\dfrac{\gamma \nu}{k+1} M \!+\! \delta A \!-\! \eta U \\[3mm]
(\eta \!+\! \delta)U \!-\! \dfrac{\gamma \nu}{k+1} M
\end{pmatrix}\!. \]

Next step is to calculate the Jacobian matrices of $\mathcal{F}, \mathcal{V}_i, i=1,2$ and evaluate them in the trivial equilibrium $\Eb_0=(0,0,0)$:

\[ F:= \left. \dfrac{\partial \mathcal{F}}{\partial \Xb} \right|_{\Eb_0} = \begin{pmatrix} 0 & 0 & r \rho \\   0 & 0 & (1-r) \rho \\ 0 & 0 & 0 \end{pmatrix}, \]
\[
V_1:= \left. \dfrac{\partial \mathcal{V}_1}{\partial \Xb} \right|_{\Eb_0} =
\begin{pmatrix}
\dfrac{\alpha k}{k+1} \!+\! \mu & 0 & 0 \\[3mm]
0 & \nu \!+\! \delta & - \eta \\[3mm]
0 & - \nu & \eta \!+\! \delta
\end{pmatrix}, \qquad  V_2:= \left. \dfrac{\partial \mathcal{V}_2}{\partial \Xb} \right|_{\Eb_0} =
\begin{pmatrix}
\dfrac{\alpha k}{k+1} \!+\! \mu & 0 & 0 \\[3mm]
\dfrac{\gamma \nu}{k+1} & \delta & - \eta \\[3mm]
- \dfrac{\gamma \nu}{k+1} & 0 & \eta \!+\! \delta
\end{pmatrix}. \]
Subsequently, the next-generation matrices of the form $F V_i^{-1}, i=1,2$ can be constructed for the dynamical systems \eqref{sys-abun} and \eqref{sys-scar}, and their  respective eigenvalues can be identified. Let us start by constructing the next-generation matrix for the system \eqref{sys-abun}:

\[ F V_1^{-1} \!=\!
\begin{pmatrix}
0 & 0 & r \rho \\[2mm]
0 & 0 & (1\!-\!r) \rho \\[2mm]
0 & 0 & 0
\end{pmatrix} \!\! \begin{pmatrix}
\dfrac{k \!+\! 1}{k(\alpha \!+\! \mu) \!+\! \mu} & 0 & 0 \\[3mm]
0 & \dfrac{\delta \!+\! \eta}{\delta (\delta \!+\! \eta \!+\! \nu)} & \dfrac{\eta}{\delta (\delta \!+\! \eta \!+\! \nu)} \\[3mm]
0 & \dfrac{\nu}{\delta (\delta \!+\! \eta \!+\! \nu)} & \dfrac{\delta \!+\! \nu}{\delta (\delta \!+\! \eta \!+\! \nu)}
\end{pmatrix} \!=\! \begin{pmatrix}
0 & \dfrac{r \rho \nu}{\delta (\delta \!+\! \eta \!+\! \nu)} & \dfrac{r \rho (\delta \!+\! \nu)}{\delta (\delta \!+\! \eta \!+\! \nu)} \\[4mm]
0 & \dfrac{(1 \!-\! r) \rho \nu}{\delta (\delta \!+\! \eta \!+\! \nu)} & \dfrac{(1 \!-\! r) \rho (\delta \!+\! \nu)}{\delta (\delta \!+\! \eta \!+\! \nu)} \\[4mm]
0 & 0 & 0
\end{pmatrix}. \]
The next-generation matrix $F V_1^{-1}$ corresponding to the dynamical system \eqref{sys-abun} is upper-triangular, and its eigenvalues are located on the main diagonal. There are two zero eigenvalues ($\lambda_1^1=\lambda_3^1=0$) and a positive one that determines the spectral radius of $F V_1^{-1}$:

\[ \lambda_2^1=\dfrac{(1 \!-\! r) \rho \nu}{\delta (\delta \!+\! \eta \!+\! \nu)} = \N_F >0.  \]

Thus, the spectral radius of the next-generation matrix $F V_1^{-1}$ corresponding to the dynamical system \eqref{sys-abun} does not depend on the external parameters $k$ and $\alpha$. Moreover, the spectral radius of $F V_1^{-1}$ is precisely the basic offspring number $\N_F$ related to the female population of psyllids, which was already derived for the dynamical system \eqref{sys-pws}-\eqref{F1} describing the natural ACP dynamics under the male abundance (cf. formula \eqref{offspring-MF}).

Furthermore, following Subsection \ref{subsec-c1}, replacing $-\mu$ by $-\mu+\alpha \dfrac{k}{k+1}$, and using similar computations to those developed in Appendix A, one can deduce that there exists a strictly positive equilibrium $\Eb_{1,P}^*\leq \Eb_{1}^*$ whose coordinates are

\begin{subequations}
\label{equilibrium-E1P}
\begin{align}
M_{1,P}^* &= \dfrac{(k+1) r \delta}{(k+1) r \delta + \big(\alpha k + (k+1) \mu \big)(1-r)} \: \dfrac{1}{\sigma} \ln \N_{F}  \\[2mm]
A_{1,P}^* &= \dfrac{(\eta+\delta )(1-r)}{ \nu+\eta+\delta\ } \: \dfrac{ \alpha k + (k+1) \mu }{ (k+1) r \delta + \big( \alpha k+ (k+1) \mu \big) (1-r) } \: \dfrac{1}{\sigma} \ln \N_{F}, \\[2mm]
U_{1,P}^* &= \dfrac{\nu(1-r)}{\nu+\eta+\delta} \: \dfrac{ \alpha k+(k+1) \mu }{ (k+1) r \delta + \big( \alpha k+(k+1)\mu \big)(1-r) } \: \dfrac{1}{\sigma} \ln \N_{F}
\end{align}
\end{subequations}
When $\alpha=0$ or $k=0$, we recover the equilibrium $\Eb_{1}^*$ given by \eqref{E-1}. Then, we derive the following result.
\begin{proposition}
\label{prop2_P}
Consider the dynamical system \eqref{sys-abun}.
\begin{itemize}
\item
Assume $\N_F<1$. Then $\Eb_0$ is LAS.
\item
Assume $\N_F>1$. Then $\Eb_{1,P}^*$ is LAS and $\Eb_0$ is unstable.
\end{itemize}
\end{proposition}
\begin{remark}
It is interesting to notice that trapping alone is insufficient to suppress the pest population drastically. When $\alpha=0$, releasing a small amount of pheromones has absolutely no impact on the population.
\end{remark}

Let us now construct the next-generation matrix for the system \eqref{sys-scar}:

\[ F V_2^{-1} \!=\!
\begin{pmatrix} 0 & 0 & r \rho \\[2mm]
0 & 0 & (1\!-\!r) \rho \\[2mm]
0 & 0 & 0
\end{pmatrix} \!\!
\begin{pmatrix}
\dfrac{k \!+\! 1}{k(\alpha \!+\! \mu) \!+\! \mu} & 0 & 0 \\[3mm]
- \dfrac{\gamma \nu}{(\delta \!+\! \eta) \big( k (\alpha \!+\! \mu) \!+\! \mu \big)} & \dfrac{1}{\delta} & \dfrac{\eta}{\delta (\delta \!+\! \eta)} \\[3mm]
\dfrac{\gamma \nu}{(\delta \!+\! \eta) \big( k (\alpha \!+\! \mu) \!+\! \mu \big)} & 0 & \dfrac{1}{\delta \!+\! \eta}
\end{pmatrix} \!=\!
\begin{pmatrix}
\dfrac{r \rho \gamma \nu}{(\delta \!+\! \eta) \big( k (\alpha \!+\! \mu) \!+\! \mu \big)} & 0 & \dfrac{r \rho}{\delta \!+\! \eta} \\[4mm]
\dfrac{(1 \!-\! r) \rho \gamma \nu}{(\delta \!+\! \eta) \big( k (\alpha \!+\! \mu) \!+\! \mu \big)} & 0 & \dfrac{(1 \!-\! r) \rho}{\delta \!+\! \eta} \\[4mm]
0 & 0 & 0 \end{pmatrix} \]
The next-generation matrix $F V_2^{-1}$ corresponding to the dynamical system \eqref{sys-scar} has only one linearly independent row (or column), and therefore it possesses only one non-zero eigenvalue:

\begin{equation}
\label{N-M-feed}
\lambda_1^2 = \dfrac{r \rho \gamma \nu}{(\delta \!+\! \eta) \big( k (\alpha \!+\! \mu) \!+\! \mu \big)}, \qquad \lambda_2^2 =\lambda_3^2 =0.
\end{equation}
Notably, this positive eigenvalue, which defines the spectral radius of the next-generation matrix $F V_2^{-1}$, depends on the external parameters $k$ and $\alpha$. On the other hand, let us recall that the spectral radius of $F V_2^{-1}$ expresses the mean number of descendants produced by one individual during its lifetime and defines the basic offspring number for the dynamical system \eqref{sys-scar}, that is,

\[ \widetilde{\N}_M (k, \alpha) := \dfrac{r \rho \gamma \nu}{(\delta \!+\! \eta) \big( k (\alpha \!+\! \mu) \!+\! \mu \big)}. \]

It is worthwhile to point out that $ \widetilde{\N}_M (k, \alpha) = \N_M$ only if $k=0$ and $\alpha=0$; otherwise, we have $\widetilde{\N}_M (k, \alpha) < \N_M,$ where $\N_M$ denotes the basic offspring number corresponding to the dynamical system \eqref{sys-pws}-\eqref{F2} describing the natural ACP dynamics under the male scarcity (cf. formula \eqref{offspring-MF}).

Let us also recall that the fulfillment of condition $\widetilde{\N}_M (k, \alpha) < 1$ would guarantee a local extinction of the ACP population described by the dynamical system \eqref{sys-scar}. Therefore, one may choose the values of parameters $k >0$ and $\alpha \in [0,1]$ to satisfy this condition, namely

\begin{equation}
\label{feed} k > k^{*} (\alpha)= \dfrac{\mu (\N_M - 1)}{\alpha + \mu},
\end{equation}
where $k^{*} (\alpha)$ is a curve decreasing with respect to $\alpha \in [0,1]$. Thus, we have established that the feedback gain $k>0$ in \eqref{gain} should satisfy the condition \eqref{feed} in order to guarantee a local extinction of the ACP population.

Similarly to the male-abundance model, we can show the existence of a positive equilibrium, $\Eb_{2,P}^*$, for the male-scarcity model \eqref{sys-scar}. Following Subsection \ref{subsec-c2}, a straightforward computation leads to

\begin{equation*}
    M_{2,P}^*+A_{2,P}^*+U_{2,P}^*=\dfrac{1}{\sigma } \ln \widetilde{\N}_M (k, \alpha).
\end{equation*}
This equality is meaningful only when $\widetilde{\N}_M (k, \alpha)>1$, that is

\begin{equation}
\label{feed_E2P} k < k^{*} (\alpha)= \dfrac{\mu (\N_M - 1)}{\alpha + \mu},
\end{equation}
which is exactly the opposite of \eqref{feed}. Further computations show that

\[ A_{2,P}^* =\left( \dfrac{(1-r) \rho}{\delta}\dfrac{(\alpha+\mu)k+\mu}{\mu \N_{M}} - 1 \right)U= \big( \widetilde{\theta}_{M}(k,\alpha) -1 \big)U
\]
exists if and only if $\widetilde{\theta}_{M}(k,\alpha) = \dfrac{(1-r) \rho}{\delta} \dfrac{1}{\widetilde{\N}_{M}(k,\alpha)} >1$. In fact, since $\theta_{M}=\widetilde{\theta}_{M}(0,0)>1$, and $\widetilde{\theta}_{M}(k,\alpha)$ being an increasing function of $k$ and $\alpha$, we deduce that $\widetilde{\theta}_{M}(k,\alpha)>1$, for all $k\geq 0, \alpha \in [0,1]$ satisfying the condition $\widetilde{\N}_M (k, \alpha) > 1$. Furthermore, the coordinates of $\Eb_{2,P}^*$ have a form similar to \eqref{equilib2} with $\theta_M$ and $\N_M$ replaced by $\widetilde{\theta}_{M}(k,\alpha)$ and $\widetilde{\N}_M (k, \alpha)$, respectively.

Finally, using the same computations as in Appendix A, we derive the following result.
\begin{proposition}
\label{prop3_P}
Consider the dynamical system \eqref{sys-scar}, and assume $\N_M>1$.
\begin{itemize}
\item
If $k^{*}(\alpha)<k$, then $\Eb_0$ is locally asymptotically stable (LAS).
\item
If  $k^{*}(\alpha)>k$ and $\theta_M > 1$, then $\Eb_{2,P}^*$ is LAS and $\Eb_0$ is unstable.
\end{itemize}
\end{proposition}
\begin{remark}
The previous result shows that elimination is reachable if the proportion of the released pheromones is sufficiently large. Indeed, emitting an amount of pheromones proportionally to the number of $A$ individuals present at each time will avoid having the Allee effect and bistability (exhibited in the open-loop case, see Subsection \ref{subsec-open}) and also allow that $\Eb_0$ be reachable and LAS even when a gradually decreasing amount of pheromones is being released.
\end{remark}

In the following section, we discuss the interplay between the choice of the strength of lure $A_p$, including the feedback gain $k$, and the male-killing rate $\alpha$ and provide illustrations of the open-loop and closed-loop approaches using numerical simulations.

\section{Numerical simulations and discussion}
\label{sec-num}

Using the parameter values from Table \ref{tab-1}, we calculate first the basic offspring numbers $\N_M$ and $\N_F$ for the natural dynamics of \textit{D. citri} (see formulas \eqref{offspring-MF}) described by the PWS dynamical system \eqref{sys-pws}-\eqref{F2}:

\[ \N_M = 37.4256, \qquad \N_F = 32.2013. \]
Thus, for the natural ACP dynamics, we have $\N_M > \N_F > 1$ meaning the abundance of males, so the population of \textit{D. citri} evolves according to the system \eqref{sys-pws}-\eqref{sys-F1}.

Let us assume that, at the initial time $t=0$, the natural ACP population is close to its steady state $\Eb_1  \in \mathbb{M}_a$, that is,

\[ M(0)= 1519 \approx M_1^{*}, \qquad A(0)=1590 \approx A_1^{*} , \qquad U(0)= 383 \approx U_1^{*}. \]

\begin{figure}[t]
\centering
	\includegraphics[width=0.9\linewidth]{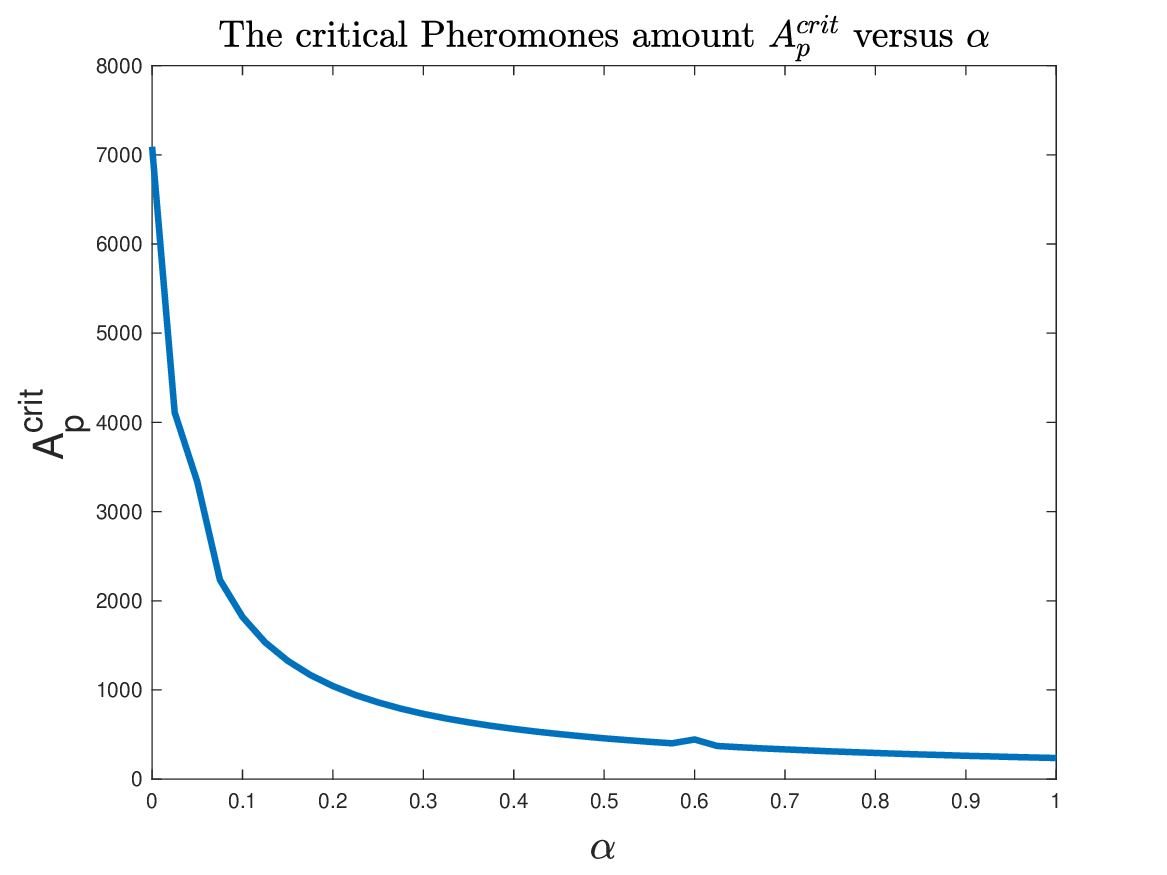}
 \caption{Open-loop control. Numerical estimates of the critical pheromones amount, $A_{p}^{crit}$, versus $\alpha$, the trapping killing rate, such that for all $A_p>A_{p}^{crit}$, the equilibrium $\Eb_0=\bf{0}$ is globally asymptotically stable (GAS) \label{fig:Ap-crit-alpha-open}
 }
\end{figure}

In the sequel, we perform numerical simulations of the PWS system \eqref{syscon} with the initial conditions given above and varying the external parameters $A_p$ and $\alpha$.

Thanks to Theorem \ref{theo_GAS_scarse_model} (page \pageref{theo_GAS_scarse_model}) formulated for the open-loop control approach and using the formula derived in Appendix C, we can estimate the critical amount of pheromones, $A_p^{crit}$, necessary to guarantee the convergence towards $\Eb_0=\bf{0}$ for different values of $\alpha$, the trapping killing rate --- see Figure \ref{fig:Ap-crit-alpha-open}. This figure shows a considerable difference between $\alpha=0$ and $\alpha>0,$ meaning that pheromones alone are not enough to control the population. Of course, the duration to enter more or less rapidly in the basin of attraction of $\Eb_0=\bf{0}$ will depend on the amount of pheromones released. Here, we can use the massive and small releases strategy developed earlier for the Sterile Insect Technique approach \cite{Anguelov2020}.

\begin{figure}[t]
\centering
	\includegraphics[width=0.9\linewidth]{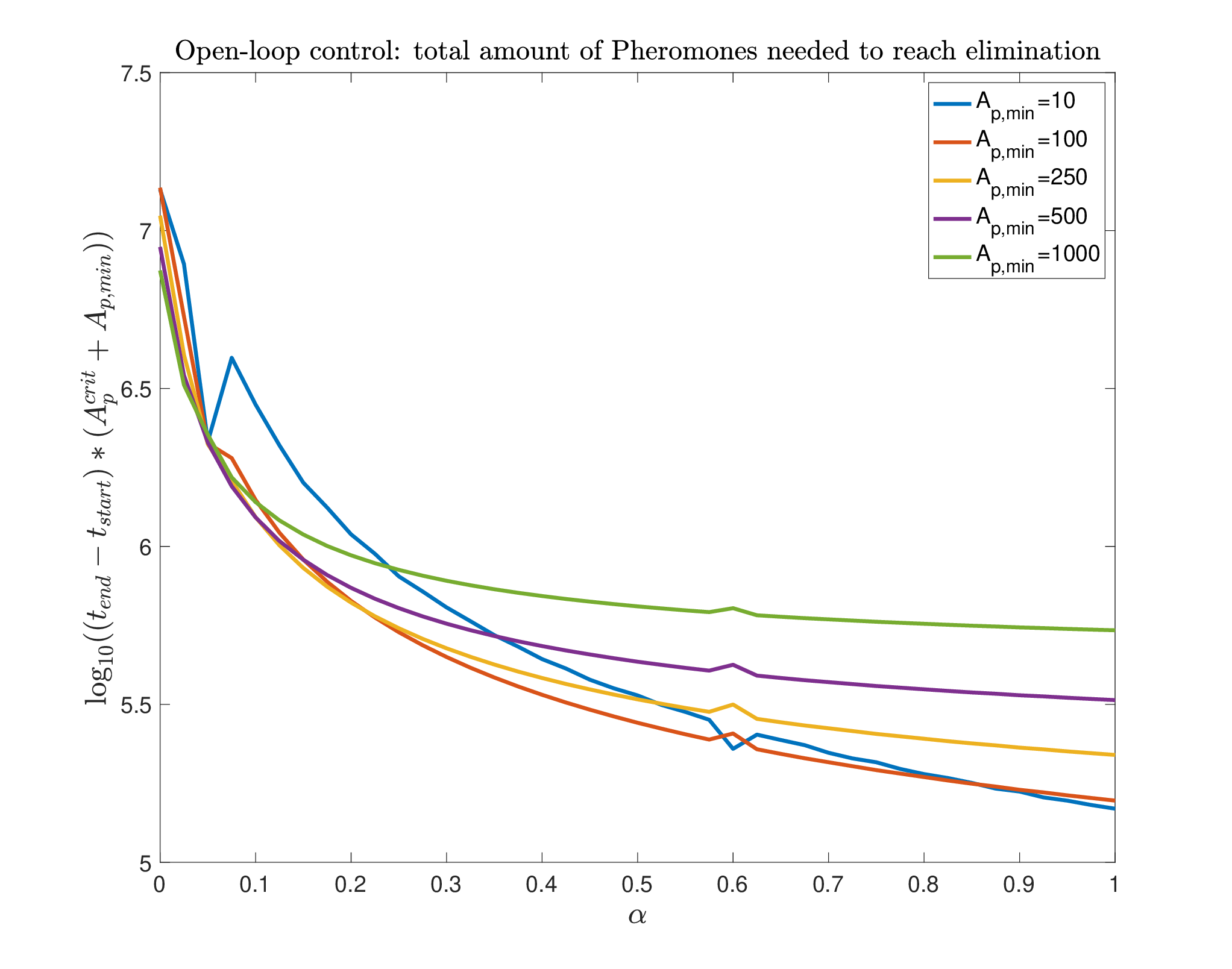}
\caption{Open-loop control. The total amount of pheromones needed to reach elimination for different values of $A_{p,min}>0$ such that $A_p=A_p^{crit}+A_{p,min}$. \label{total-open-loop-control}
}
\end{figure}

In Figure \ref{total-open-loop-control}, we show the total amount of pheromones needed to reach elimination under the open-loop control approach where $A_p$ is chosen above $A_p^{crit}$, i.e., $A_p=A_p^{crit}+A_{p,min}$, for different values of $A_{p,min}$. As expected, the lower $A_{p,min}>0$, the lowest the total pheromones amount, but the intervention becomes longer. Thus, there is a balance to be found between the amount of pheromones available for releasing and the duration of the treatment.

\begin{figure}[t]
\centering
	\includegraphics[width=0.9\linewidth]{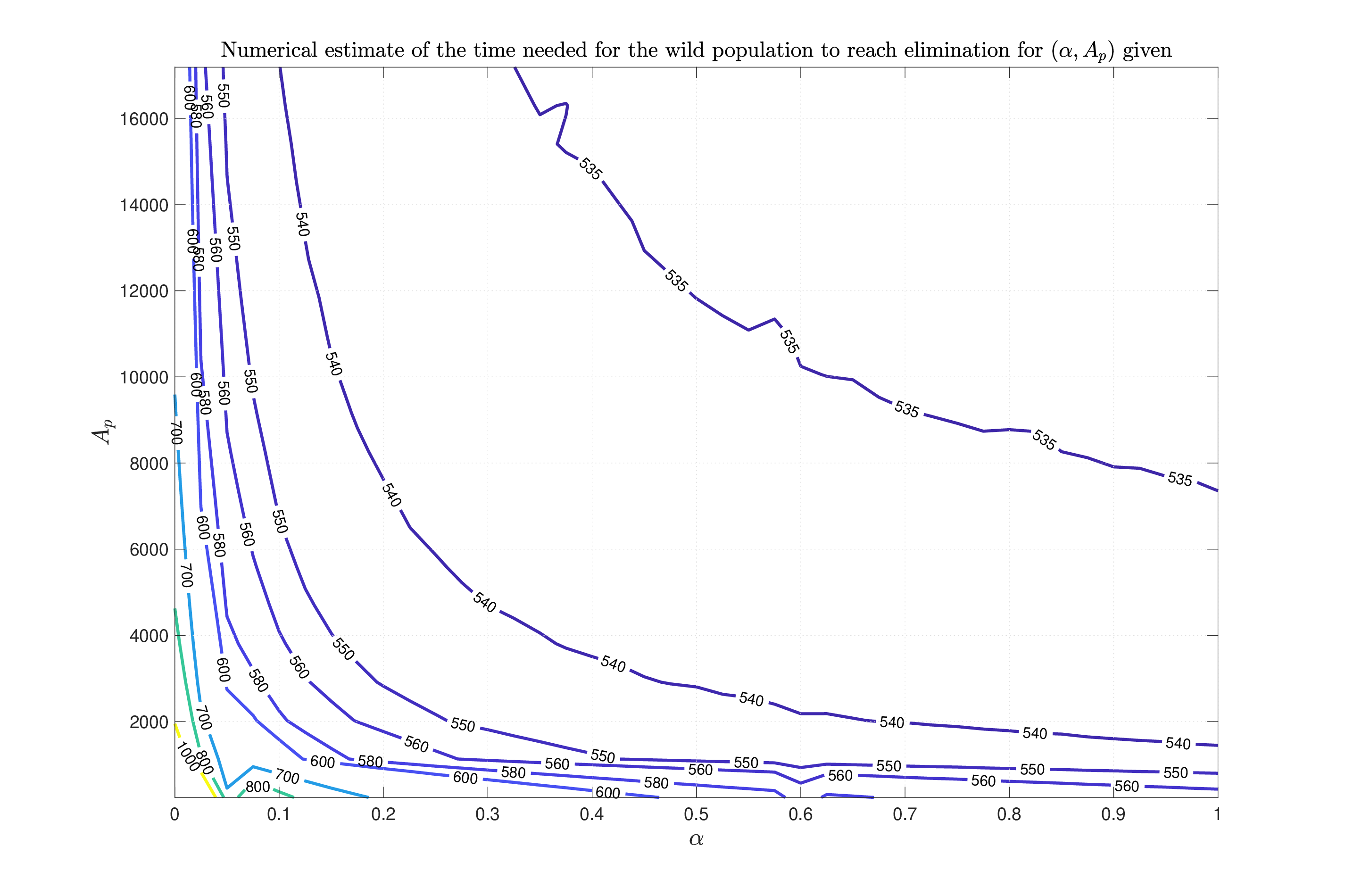}
 \caption{Open-loop control. Minimal time estimates for  (nearly) eliminating the ACP population, initially at equilibrium, with $\big\|(M,A,U)^T \big\|_{\infty}<10^{-1}$, induced by the trapping rate $\alpha$ and the amount  $A_p$ of released pheromones. \label{fig:time}
 }
\end{figure}

Then,  Figure \ref{fig:time} shows that combining pheromones and trapping is also essential to lower the time needed to (nearly) reach elimination. In fact, the lowest time value is $t=535$ days even if $A_p$ is very large and $\alpha$ is close to $1$. This is an interesting result because it shows that even if a sufficiently large value for $A_p$ is available (to have $\Eb_0= \bf{0}$ as a global attractor), using it in vast quantities will not be helpful. Of course, the larger $\alpha$, the smaller $A_p$: this shows that there is a tradeoff between these two quantities. In any case, it is deduced that releasing the pheromones alone is impractical.

\begin{figure}[t]
\centering
	\includegraphics[width=0.6\linewidth]{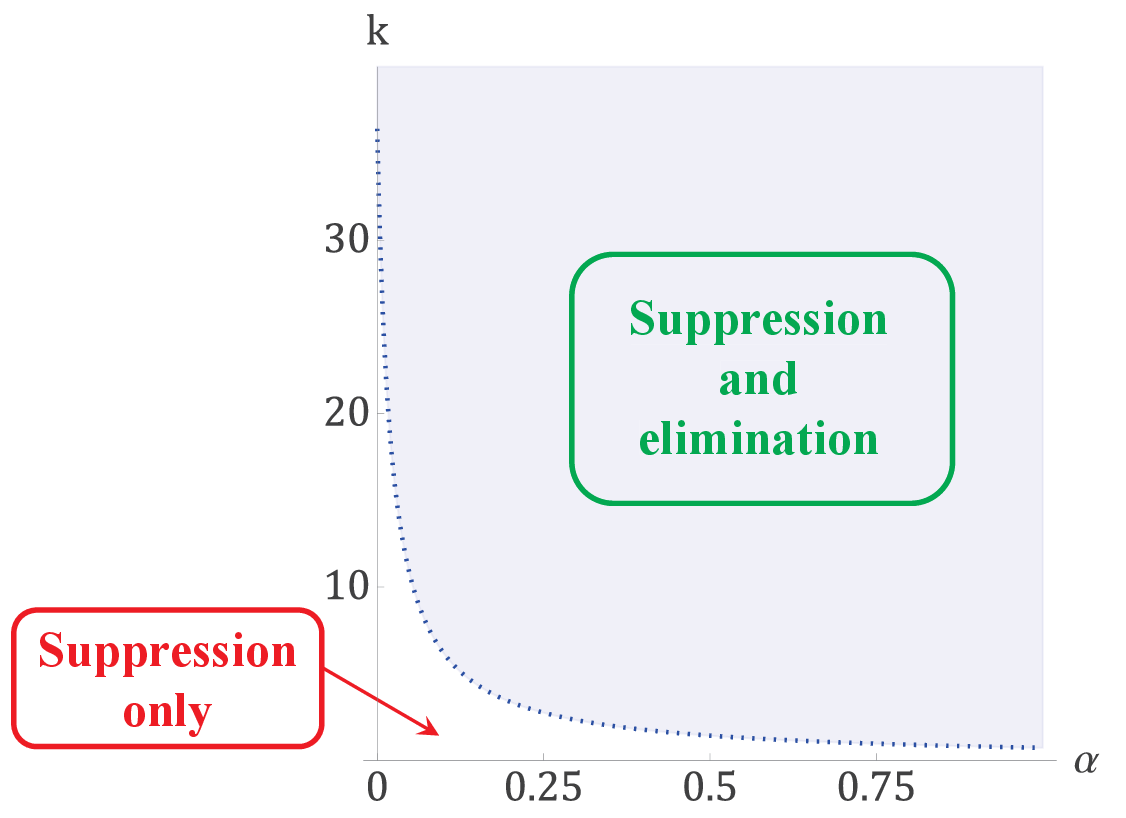}
\caption{Closed-loop control. The choice of the feedback gain $k$ and male-killing rate $\alpha$ according to the condition \eqref{feed} \label{fig-k}}
\end{figure}

Turning to the closed-loop control approach, Figure \ref{fig-k} presents the plot of the curve $k^{*} (\alpha)$ (blue dotted line) that divides the positive quadrant into two regions. The unshaded region below the curve $k^{*} (\alpha)$ contains all the values of $k$ and $\alpha$ that guarantee only the suppression of the local ACP population. In contrast, the shaded region above the curve $k^{*} (\alpha)$ contains all the values of $k$ and $\alpha$ that guarantee the elimination of the local ACP population.

\begin{figure}[t]
\centering
\begin{tabular}{cc}
\hspace{-20mm} $k < k^*(\alpha)$ &  $k > k^*(\alpha)$ \\
& \\
\includegraphics[height=5cm]{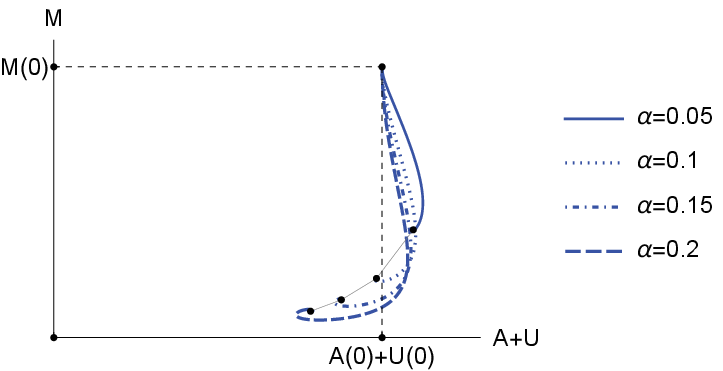} & \includegraphics[height=5cm]{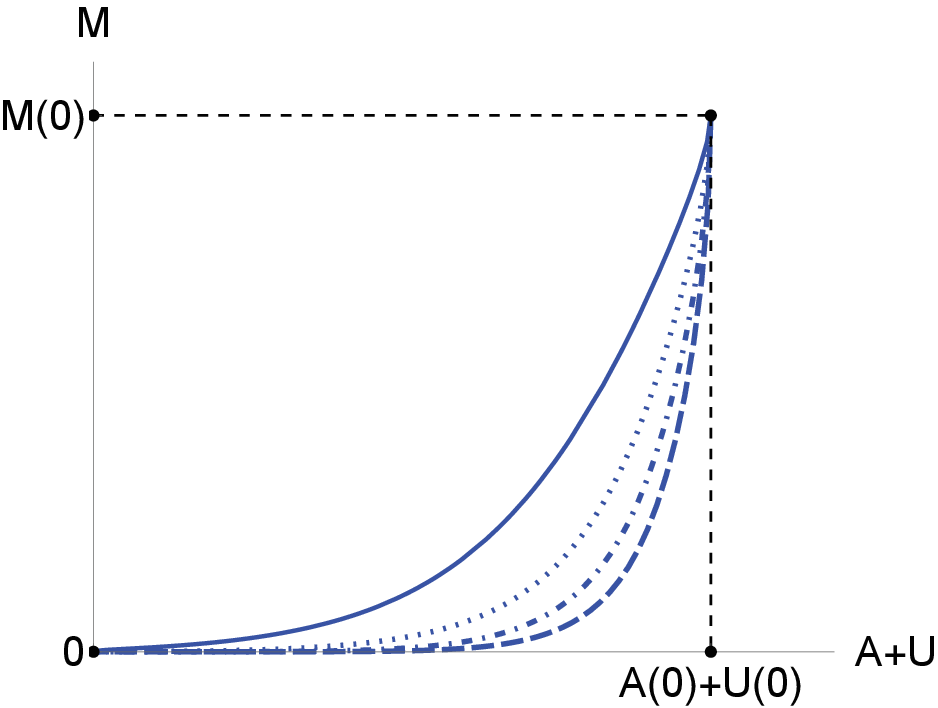} \\
\end{tabular}
\caption{Closed-loop control. The effect of male-killing rate $\alpha$ on the phase population trajectories of males $M$ and females $F=A+U$ in the case of population suppression ($k < k^*(\alpha)$, \textit{left chart}) and local elimination ($k > k^*(\alpha)$, \textit{right chart})  \label{fig:phase}}
\end{figure}

From Figure \ref{fig-k}, one may deduce that local elimination of the ACP population is possible even if $\alpha =0$ (meaning that none of the male insects are killed when approaching or entering the trap). However, in such a case, the feedback gain $k$ should be substantial ($k \geq 36.43$), meaning that a vast amount of sex pheromone should be emitted ($A_p \geq 36.43 A(t)$, in terms of the ``false'' females). On the other hand, Figure \ref{fig-k} also displays that the effect of male-killing rate $\alpha$ on the reduction of the total amount of sex pheromone needed is more visible for smaller values of $\alpha$ (below $50\%$) than for its larger values (above $50\%$). Nonetheless, even the traps with $100\%$ male-killing rate ($\alpha =1$) will still need a small feedback gain of about $k=1$ ($A_p \geq A(t)$, in terms of the ``false'' females) to reach an eventual local elimination of \textit{D. citri}.

Figure \ref{fig:phase} shows the evolution of the male $M$ and female $F=A+U$ population in the phase plane $(M, F)$ for different values of the male-killing rate $\alpha$. When the feedback gain is below the curve $k^{*}(\alpha)$ (here we have chosen $k=2.5$), we observe the convergence of the phase trajectories to another positive equilibrium $\Eb_{1, P}^{*}$ whose coordinates are given by \eqref{equilibrium-E1P}. For different values of $\alpha$, the corresponding $\Eb_{1, P}^{*}$ are marked by the black points on the left chart of Figure \ref{fig:phase}. Notably, this equilibrium moves closer to the origin ($\Eb_0$) as $\alpha$ increases, and the population of males decays faster than that of females.

On the other hand, if the feedback gain is above the curve $k^{*}(\alpha)$ (here we have chosen $k=\N_M$), we observe the convergence of the phase trajectories to the trivial equilibrium $\Eb_0$ (see the right chart in Figure \ref{fig:phase}). Here, the density of males also decays faster than the density of females.

\begin{figure}[t]
\centering
	\includegraphics[width=0.95\linewidth]{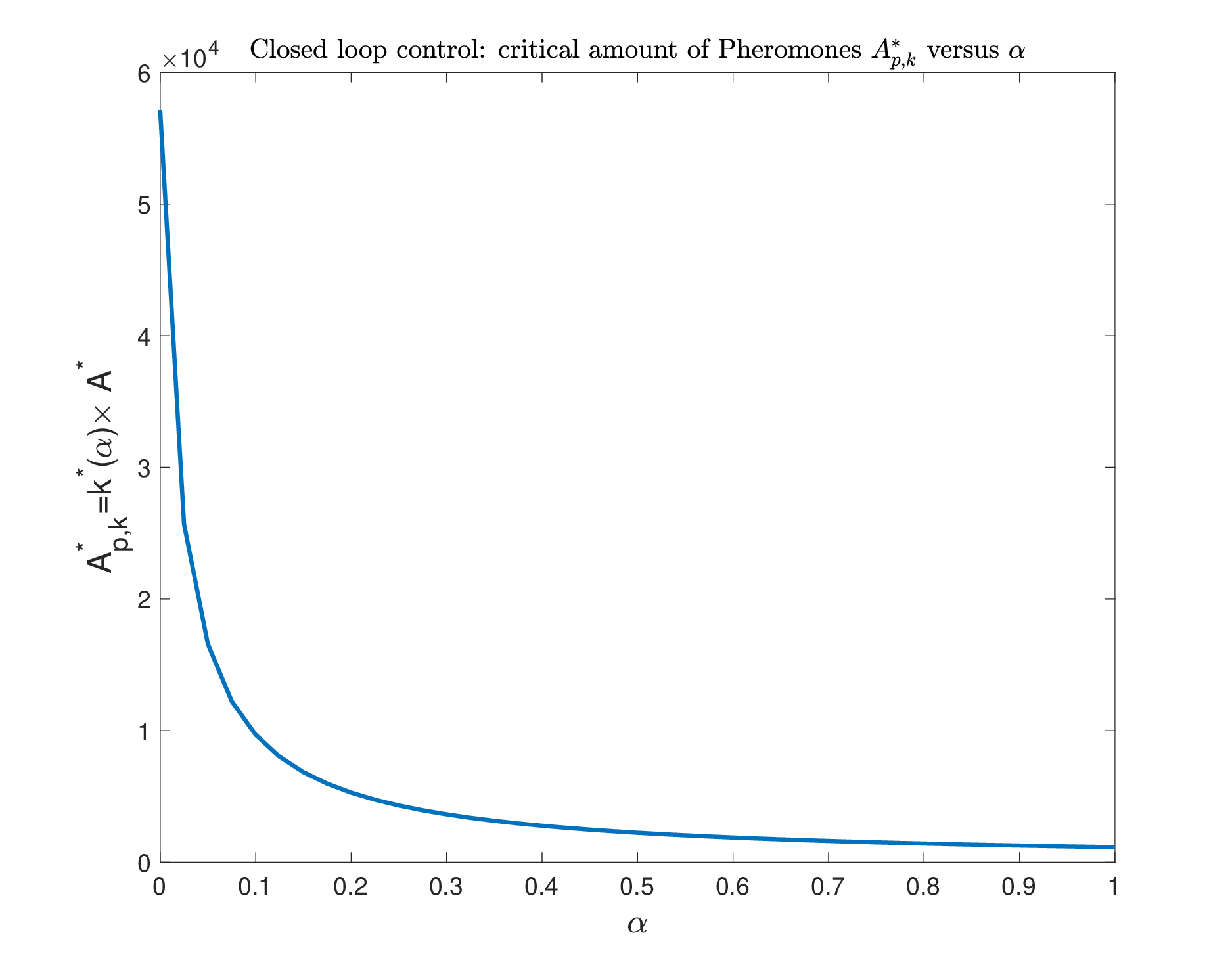}
 \caption{Closed-loop control. Estimate of the critical pheromones amount, $A_{p,k}^*=k^*(\alpha) \times A^*$, versus $\alpha$, the trapping killing rate, such that for all $A_p>A_{p,k}^*$, the equilibrium $\Eb_0=\bf{0}$ is globally asymptotically stable \label{Ap_crit_versus_alpha_closed_loop}}
\end{figure}

To compare the closed-loop and open-loop approaches, at least at the beginning of the pheromones treatment, we show the amount of pheromones needed to start the closed-loop control: see Figure \ref{Ap_crit_versus_alpha_closed_loop} (page \pageref{Ap_crit_versus_alpha_closed_loop}). Thus, contrasting this figure with Figure \ref{fig:Ap-crit-alpha-open} (page \pageref{fig:Ap-crit-alpha-open}), it is straightforward to see that the closed-loop control is very costly compared to the open-loop control. Indeed, the closed-loop control requires almost $9$ times more pheromones than the open-loop control, regardless of the values for $\alpha$.

Notably, the result presented in Figure \ref{Ap_crit_versus_alpha_closed_loop} is somewhat idealistic because it corresponds to ``continuous estimations'' of the population size of female psyllids $A(t)$, available for mating, assuming that for all $t \geq 0$ the size of $A(t)$ can be accurately assessed. In practice, however, estimating the size of a local insect population can be a challenging and expensive task. Therefore, continuous population size estimations are unfeasible, and most Integrated Pest Management (IPM) programs conduct such assessments with different frequencies but not more often than every two weeks.

\begin{figure}[t]
\centering
	\includegraphics[width=0.95\linewidth]{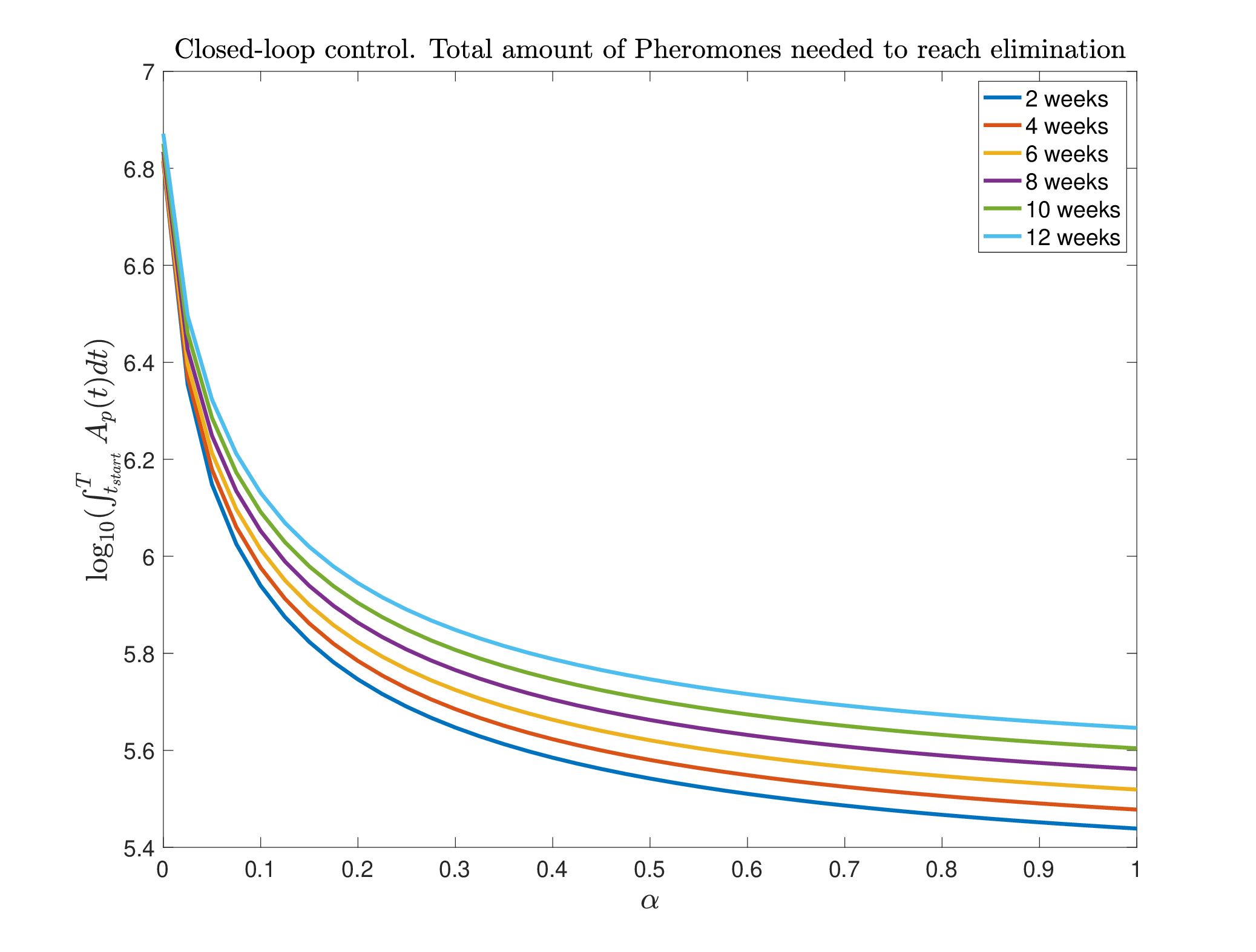}
\caption{Closed-loop control. Total amount of pheromones needed to reach elimination, when the population is estimated every $2n$ weeks, where $n=1,2...,6$. \label{total-closed-loop-control}}
\end{figure}

Using the closed-loop control approach, we have assessed the total amount of pheromones $A_p$ (also for different values of the male-killing rate $\alpha$) needed to reach elimination if the size of female psyllids available for mating, $A(t)$, is estimated every  $2n$ weeks, where $n=1, 2, \ldots, 6.$ The underlying results are displayed in Figure \ref{total-closed-loop-control}, where $A_p= \big( k^*(\alpha)+ 1 \big)\times A(t_{j}^{2n})$, and $t_{j}^{2n}$ with $j=1,2,...$, denote the times when the wild population is estimated.

According to definition \eqref{gain}, closed-loop control becomes useful once the population is or has become small enough, meaning that $A_p = kA(t)$ is not too large. This rationale leads us to consider a mixed-type control and choose $A_p$, for instance, in the following way:

\[ A_p = \min \Big\{ A_{p}^{crit}(\alpha)+A_{p,min}; \big( k^{*}(\alpha)+1 \big) A(t^*) \Big\} \]
with $A_{p,min}=500$, for a given $\alpha$, where $A(t^*)$ is an estimated value of $A$ at a given time $t^*$. Like in the closed-loop case, we estimate the population size every $2n$ weeks, where $n=1, 2, \ldots, 6$. Thus we choose

\[ A_p(t_{j}^{2n})=\min \Big\{  A_{p}^{crit}(\alpha) + 500; \big( k^{*}(\alpha) + 1 \big) A(t_{j}^{2n}) \Big\}, \qquad j=1,2,... \]
and $n=1,2,..,6$.

\begin{figure}[t]
\centering
	\includegraphics[width=0.95\linewidth]{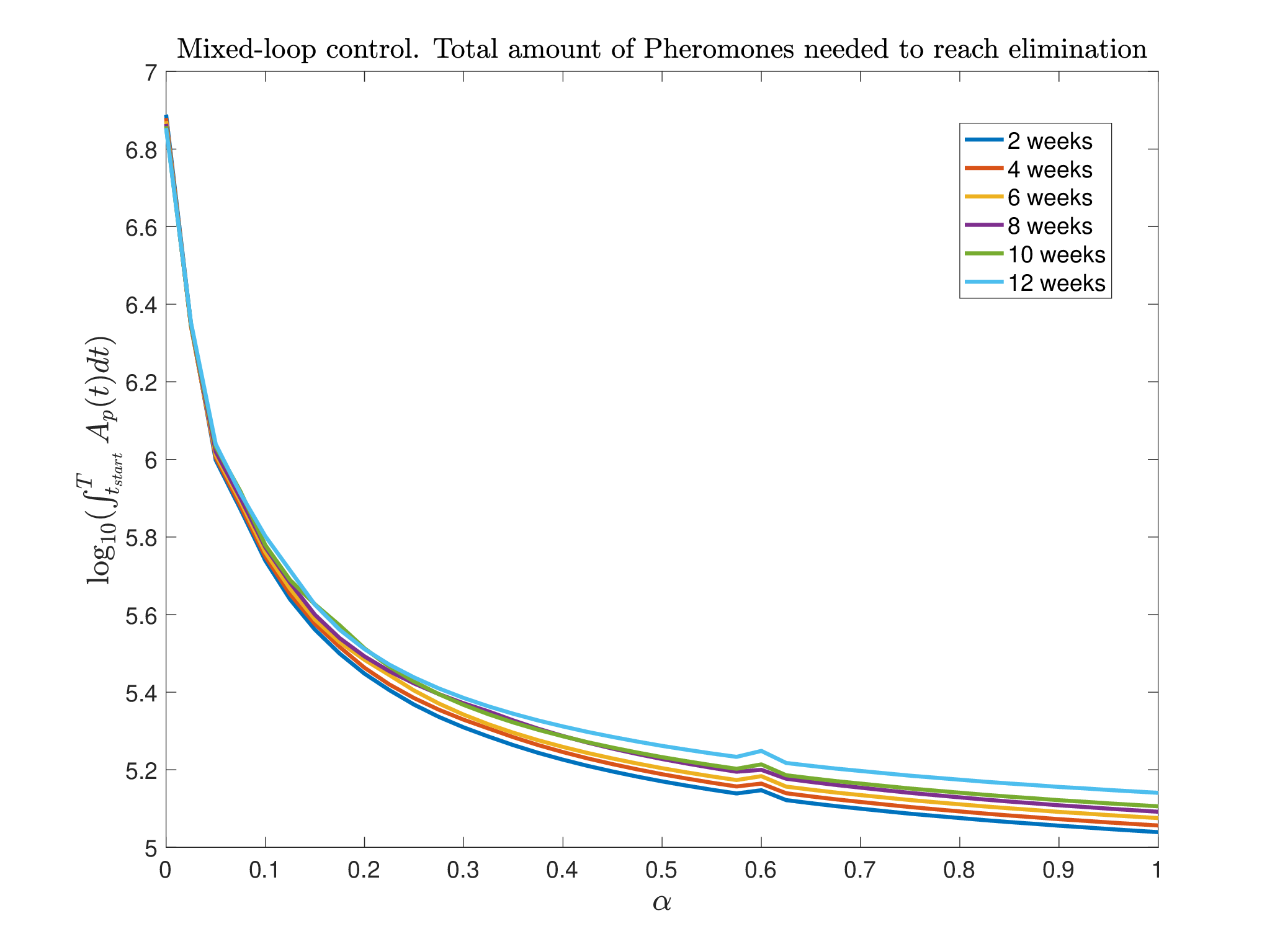}
\caption{Mixed-type control. Total amount of pheromones needed to reach elimination under an ``open-loop -- closed-loop'' control approach, when the population is estimated every $2n$ weeks, where $n=1,2...,6$. \label{mixed-control}}
\end{figure}

To illustrate this approach, we provide in Figure \ref{mixed-control} (page \pageref{mixed-control}) the total amount of pheromones, for different values of $\alpha$, related to the estimates of the population size carried out every $2n$ weeks, where $n=1,2,...,6$. As expected, the mixed-type control is functional when $\alpha$ is large enough (compare with Figure \ref{total-open-loop-control}). Last, the smaller is $n$, the less is the total amount of pheromones for any $\alpha$. However, we must be aware that estimating the population size in the field can be very difficult. Altogether, the mixed-type control provides the best result.

\section{Conclusions}
\label{sec-concl}

Citrus fruits are among the most important crops worldwide, and many citrus cultures worldwide face the threat of Huanglongbing (HLB) or citrus greening disease \cite{Bove2006}. This disease is mainly transmitted by the Asian citrus psyllid, \textit{Diaphorina citri}, an invasive psyllid species that colonizes citrus orchards in different parts of the world \cite{Aidoo2022}. Controlling this pest population is a challenging task, and Integrated Pest Management (IPM) programs are seeking environmentally friendly strategies that may replace the traditional ones based on pesticides. In this context, using pheromone traps seems rather promising because the attraction and direct removal of male insects induce mating disruption, thus reducing future offspring and suppressing the overall pest population.

In this paper, we proposed and analyzed a model formulated as a piecewise smooth ODE system that describes the natural population dynamics of Asian citrus psyllids. This model was further amended with the external control actions, expressed by the parameters $A_p$ (the strength of lure) and $\alpha$ (traps male-killing rate), to mimic the introduction of sex pheromone traps that enforce mating disruption and may lead either to the suppression or local extinction of the ACP population. From the theoretical standpoint, the choice of external parameters $A_p$ and $\alpha$, as well as their interplay, was conducted based on two operational control modes, the open-loop control approach (operating on a predefined sequence of actions) and the closed-loop control approach (whose actions can respond to changes in the system behavior). For both techniques, we have identified a critical curve or mapping, that is, $A_p$ as a function of $\alpha$, that plays the role of a threshold, below which only the ACP population suppression is achievable and above which the local pest extinction can be attained. These theoretical findings also allowed us to perform  qualitative \textit{in silico} testings of the model and estimate not only the total amount of sex pheromones needed for the local elimination of the ACP population but also the minimum time to reach a local extinction of this pest.

Summarizing the outcomes of the present work, we would like to highlight the following insight regarding the control of \textit{D. citri} using the pheromone traps:
\begin{itemize}
\item
Pheromone traps are a reliable alternative to pesticide use because they can suppress and even eliminate the ACP pest population while producing no harm to the crop.
\item
Releasing the pheromones alone without male-killing  (i.e.,  with $\alpha=0$) may also be employed to reach the final goal of pest population suppression or elimination. Nonetheless, it would require emitting a massive quantity of pheromones, while using male-killing traps would considerably reduce the pheromone quantity needed for controlling the pest population.
\item
Increasing the male-killing rate $\alpha$ of the traps may reduce the overall costs (i.e., the total amount of pheromones needed for intervention) and the time to reach elimination.
\item
Open-loop control strategies show better results when applied to established pest populations bearing large sizes. In contrast, closed-loop control strategies perform better for emerging populations or when the population is small. Thus, combining these two control approaches renders the best overall results and requires a smaller amount of pheromones.
\end{itemize}

It is worthwhile to point out that the model proposed in this paper lays solid grounds for combining the pheromone traps with other control intervention measures as a part of the IPM programs. The latter can be tested to enhance the control efficiency against Huanglongbing, using a limited amount of pheromones, possibly by designing the optimal control strategies. Furthermore, extensions of the present work may include the transmission of infection caused by \textit{Candidatus} Liberibacter spp. and lead to the formulation and study of epidemiological models to reduce the risk of HLB spreading and its damage to citrus crops.

Finally, we hope that the outcomes of this study will provide valuable insights for developing alternatives to pesticides and also shed some light on the practical implementation of ecologically friendly pest management conducted in field trials.

\section*{Acknowledgements}

This research was funded by the National Fund for Science, Technology, and Innovation (Autonomous Heritage Fund \textit{Francisco Jos\'e de Caldas}) by way of the Research Program No. 1106-852-69523, Contract: CT FP 80740-439-2020 (Colombian Ministry of Science, Technology, and Innovation --- Minciencias), Grant IDs: CI-71242 (Universidad del Valle), 20 INTER 356 P2 (Universidad Autonoma de Occidente).
Yves Dumont is (partially) supported by the DST/NRF SARChI Chair in Mathematical Models and Methods in Biosciences and Bioengineering at the University of Pretoria, South Africa (Grant 82770). Yves Dumont acknowledges the support of the {\it Conseil R\'egional de la R\'eunion} (France), the {\it Conseil D\'epartemental de la R\'eunion} (France), the European Agricultural Fund for Rural Development (EAFRD), the European Regional Development Fund (ERDF), and the \textit{Centre de Coop\'eration Internationale en Recherche Agronomique pour le D\'eveloppement} (CIRAD), France.

\bibliographystyle{plain}
\addcontentsline{toc}{section}{References}
\bibliography{Diaphorina_CDV_Oct_2023}

\begin{appendix}
\appendix
\renewcommand{\theequation}{A-\arabic{equation}}
  \setcounter{equation}{0}  
\renewcommand{\thefigure}{A.\arabic{figure}}
  \setcounter{figure}{0}  
  \renewcommand{\thetable}{A.\arabic{table}}
   \setcounter{table}{0}  

\section*{Appendix A: proof of Propositions \ref{prop2} and \ref{prop3}}
\addcontentsline{toc}{section}{Appendix A: proof of Proposition \ref{prop2} and \ref{prop3}}
\label{appendixA}
\textbf{Proof of Proposition \ref{prop2}.} First, we compute the Jacobian related to system \eqref{sys-F1} with $\Phib_1$ defined by \eqref{F1}:

\begin{align*}
J(\Xb) & = \: \begin{pmatrix}
J_{11} & J_{12} & J_{13} \\ J_{21} & J_{22} & J_{23} \\ J_{31} & J_{32} & J_{33}
\end{pmatrix} \\
& \\
& = \: \begin{pmatrix}
-\mu \!-\! r\rho \sigma U e^{-\sigma(M+A+U)}   & - r \rho \sigma U e^{-\sigma(M+A+U)}                       & r \rho (1 \!-\! \sigma U) e^{-\sigma(M+A+U)} \\
& & \\
-(1 \!-\! r) \rho \sigma  U e^{-\sigma(M+A+U)} & - (\nu \!+\! \delta ) \!-\! (1 \!-\! r) \rho \sigma U e^{-\sigma(M+A+U)}  & \eta \!+\!  (1 \!-\!r)\rho (1 \!-\! \sigma U) e^{-\sigma(M+A+U)} \\
& & \\
0 & \nu & - (\eta \!+\!\delta )
\end{pmatrix},
\end{align*}
where $\Xb=(M,A,U)$. Thus, it is easy to compute that

\[ J(\Eb_0)=\left(\begin{array}{ccc}
-\mu & 0 & r\rho\\
0 & -\left(\nu+\delta\right) & \eta+\left(1-r\right)\rho\\
0 & \nu & -\left(\eta+\delta\right)
\end{array}\right), \]
and to show that the characteristic polynomial is given by

\begin{align*}
P_1^0 ( \lambda) & = \:  - (\mu + \lambda )\Big[ \lambda^2 + (\eta + \nu + 2 \delta) \lambda + (\eta + \delta) (\nu + \delta) - \big( \eta+ (1-r) \rho \big) \nu  \Big] \\
& = \: - (\mu + \lambda )\Big[ \lambda^2 + ( \nu + \eta +2 \delta) \lambda + \delta(\nu + \eta + \delta) - (1-r) \rho \nu \Big] \\
& = \: - (\mu + \lambda )\Big[ \lambda^2 + (\nu + \eta + 2 \delta) \lambda + \delta(\nu + \eta + \delta) \big( 1- \N_F \big) \Big]
\end{align*}
When $\N_F<1$, this polynomial  has three roots with negative real parts meaning that $\Eb_0$ is LAS. Alternatively,  $p_1(\lambda)$ has one root with positive real part when $\N_F>1$ meaning that $\Eb_0$ is a saddle point (not a repeller). In effect, a trajectory engendered by the initial condition $M(0) >0, A(0)=0, U(0)=0$ converges to $\Eb_0$ even if $\N_F>1$.

Let us now show that $\Eb_1^*$ is LAS when $\N_F>1$. The Jacobian evaluated at $\Eb_1^*$ is

\begin{equation}
\label{A-JE1}
 J \big( \Eb_1^{*} \big)=\begin{pmatrix}
-\mu - \dfrac{r \rho }{\N_F} \sigma U_1^*   & - \dfrac{r \rho }{\N_F} \sigma U_1^*                   & \dfrac{r \rho}{\N_F} \big( 1-\sigma U_1^* \big) \\
& & \\
- \dfrac{(1-r) \rho }{\N_F} \sigma U_1^*    & -(\nu+\delta) - \dfrac{(1-r) \rho }{\N_F} \sigma U_1^* & \eta + \dfrac{(1-r)\rho}{\N_F} \big( 1-\sigma U_1^* \big) \\
& & \\
0 & \nu & - (\eta+\delta )
\end{pmatrix}
\end{equation}
The characteristic polynomial of $J \big(\Eb_1^* \big)$ has the form

\[ P_1^* (\lambda) =  \lambda^3 + a_1 \lambda^2 + a_2 \lambda + a_3 ,  \]
where $a_1 = - \text{\texttt{Tr} } J \big(\Eb_1^{*} \big), a_3 = - \det J \big(\Eb_1^{*} \big),$ and

\begin{equation}
 \label{A-a2}
a_2 = \det \begin{vmatrix} J_{11} & J_{12} \\ J_{21} & J_{22} \end{vmatrix} + \det \begin{vmatrix} J_{11} & J_{13} \\ J_{31} & J_{33} \end{vmatrix} + \det \begin{vmatrix} J_{22} & J_{23} \\ J_{32} & J_{33} \end{vmatrix} := \Delta_1^{(1)} + \Delta_2^{(1)} + \Delta_3^{(1)}
\end{equation}
with $J_{ij}, i,j=1,2,3$ denoting the elements of \eqref{A-JE1}. According to Routh-Hurwitz criterion (see, e.g., \cite{Murray2002}), all roots of $P_1^*$ have negative real parts if and only if the following conditions are satisfied:

\begin{equation}
\label{A-RHsys1}
 a_1 >0, \quad a_3 > 0,   \quad \text{and} \quad a_1 a_2 - a_3 >0.
\end{equation}
Let us now check these conditions. First we note that

\[ \text{\texttt{Tr} } J \big(\Eb_1^{*} \big) = - \mu - \dfrac{r \rho}{\N_F} \sigma U^*_1 - (\nu + \delta) - \frac{(1-r) \rho}{\N_F}  \sigma U^*_1 - (\eta + \delta) < 0, \]
and thus we have

\begin{equation}
\label{A-a1}
a_1 = - \text{\texttt{Tr} } J \big(\Eb_1^{*} \big) = \mu  + \nu + \eta + 2\delta +\dfrac{\rho}{\N_F} \sigma U^*_1 > 0.
\end{equation}

To compute $\det J\big( \Eb_1^{*} \big)$, we observe that

\[  \det J \big( \Eb_1^{*} \big) = - \nu \Delta_{32}^{(1)} - (\eta+\delta ) \Delta_{33}^{(1)}, \]
where $\Delta_{32}^{(1)}, \Delta_{33}^{(1)}$ are minors of $J \big( \Eb_1^{*} \big)$ obtained by elimination of the third row and either second or third column from \eqref{A-JE1}. Effectively, we have

\begin{align*}
\Delta_{32}^{(1)} & = - \left( \mu + \dfrac{r \rho }{\N_F} \sigma U_1^* \right) \times \left[ \eta + \dfrac{(1-r)\rho}{\N_F} \big( 1-\sigma U_1^* \big) \right] + \dfrac{(1-r) \rho }{\N_F} \sigma U_1^* \times \dfrac{r \rho}{\N_F} \big( 1-\sigma U_1^* \big) \\
& = \mu \dfrac{(1-r) \rho  }{\N_F} \sigma U_1^* - \eta \dfrac{r \rho  }{\N_F} \sigma U_1^* - \mu \dfrac{(1-r) \rho  }{\N_F} - \mu \eta, \\
\Delta_{33}^{(1)} & = \left( \mu + \dfrac{r \rho }{\N_F} \sigma U_1^* \right) \times \left[ (\nu+\delta) + \dfrac{(1-r) \rho }{\N_F} \sigma U_1^* \right] - \dfrac{r \rho }{\N_F} \sigma U_1^* \times \dfrac{(1-r) \rho }{\N_F} \sigma U_1^* \\
& =  \mu \dfrac{(1-r) \rho }{\N_F} \sigma U_1^* +(\nu+\delta) \dfrac{r \rho  }{\N_F} \sigma U_1^*  + \mu (\nu+\delta).
\end{align*}
Then using the relationship

\begin{equation}
\label{A-rel1}
(\eta+\delta )(\nu+\delta) - \eta \nu = \delta (\delta + \eta + \nu) = \dfrac{(1-r) \rho \nu }{\N_F},
\end{equation}
we obtain

\begin{align*}
\det J \big( \Eb_1^{*} \big) & = - \nu \left[ \mu \dfrac{(1-r) \rho  }{\N_F} \sigma U_1^* - \eta \dfrac{r \rho  }{\N_F} \sigma U_1^* - \mu \dfrac{(1-r) \rho}{\N_F} - \mu \eta \right] \\
& \;\;\;\; -  (\eta+\delta ) \left[ \mu \dfrac{(1-r) \rho }{\N_F} \sigma U_1^* +(\nu+\delta) \dfrac{r \rho  }{\N_F} \sigma U_1^*  + \mu (\nu+\delta) \right] \\
& = -\mu (\nu + \eta + \delta)  \dfrac{(1-r) \rho }{\N_F} \sigma U_1^* - \delta (\nu + \eta + \delta) \dfrac{r \rho  }{\N_F} \sigma U_1^* + \mu  \dfrac{(1-r) \rho \nu }{\N_F} -  \mu  \delta (\nu + \eta + \delta) \\
& = - \Big[ (1-r)\mu + r \delta \Big] (\nu + \eta + \delta) \dfrac{\rho }{\N_F} \sigma U_1^*=  -  (\nu + \eta + \delta) \dfrac{\rho \vartheta}{\N_F} \sigma U_1^* < 0.
\end{align*}
Thus, we have

\begin{equation}
\label{A-a3}
a_3 = - \det J \big( \Eb_1^{*} \big) = \dfrac{\rho \vartheta}{\N_F} (\nu + \eta + \delta ) \sigma U_1^* > 0.
\end{equation}

To compute the coefficient $a_2$, we evaluate $\Delta_i^{(1)}, i=1,2,3$ that appear in \eqref{A-a2}:

\begin{align*}
\Delta_1^{(1)} & = \left( \mu + \dfrac{r \rho}{\N_F} \sigma U_1^{*} \right) \times \left( \nu + \delta + \dfrac{(1-r) \rho}{\N_F} \sigma U_1^{*} \right) - \dfrac{r \rho}{\N_F} \sigma U_1^{*} \times \dfrac{(1-r) \rho}{\N_F} \sigma U_1^{*} \\
& = \mu (\nu + \delta) + \mu \dfrac{(1-r) \rho}{\N_F} \sigma U_1^{*} + (\nu + \delta)\dfrac{r \rho}{\N_F} \sigma U_1^{*}, \\
\Delta_2^{(1)} & = (\eta + \delta) \left( \mu + \dfrac{r \rho}{\N_F} \sigma U_1^{*} \right), \\
\Delta_3^{(1)} & = (\eta + \delta) \left( \nu + \delta + \dfrac{(1-r) \rho}{\N_F} \sigma U_1^{*} \right) -\nu \left( \eta + \dfrac{(1-r) \rho}{\N_F} (1- \sigma U_1^{*}) \right)
\end{align*}
Then according to \eqref{A-a2} and using \eqref{A-rel1} we have

\begin{align*}
a_2 & = \mu (\nu + \delta) + \mu (\eta + \delta) + (\eta+\delta )(\nu+\delta) - \eta \nu - \dfrac{(1-r) \rho \nu }{\N_F} \\
& \hspace{4mm} + \dfrac{\rho }{\N_F} \sigma U_1^*  \Big[  (1-r)\mu + r (\nu + \delta) + r (\eta + \delta) + (1-r)(\eta + \delta) + (1-r) \nu \Big] \\
& = \mu (\nu + \delta) + \mu (\eta + \delta) + (\vartheta + \nu + \eta + \delta) \dfrac{\rho}{\N_F} \sigma U_1^* >0.
\end{align*}

Our goal now is to show that $a_1 a_2 - a_3 >0$. Before proceeding, we rewrite $a_1$ and $a_2$ in terms of $a_3$ using the formula \eqref{A-a3}:

\[ a_1 = (\nu + \eta + \delta) + (\mu +\delta) + \dfrac{a_3}{\vartheta (\nu + \eta + \delta)}, \quad a_2 = \mu (\nu + \eta + 2\delta) + \dfrac{a_3}{\nu + \eta + \delta} + \dfrac{a_3}{\vartheta}, \]
so that

\begin{align*}
a_1 a_2 - a_3 & = \left( (\nu + \eta + \delta) + (\mu +\delta) + \dfrac{a_3}{\vartheta (\nu + \eta + \delta)} \right) \left( \mu (\nu + \eta + 2\delta) + \dfrac{a_3}{\vartheta} + \dfrac{a_3}{\nu + \eta + \delta} \right) - a_3 \\
& = (\nu + \eta + \delta) \left( \mu (\nu + \eta + 2\delta) + \dfrac{a_3}{\vartheta} \right) + a_3 \\
& \hspace{4mm} + \left( (\mu +\delta) + \dfrac{a_3}{\vartheta (\nu + \eta + \delta)} \right) \left( \mu (\nu + \eta + 2\delta) + \dfrac{a_3}{\vartheta} + \dfrac{a_3}{\nu + \eta + \delta} \right) - a_3 \\
& = (\nu + \eta + \delta) \left( \mu (\nu + \eta + 2\delta) + \dfrac{a_3}{\vartheta} \right)  \\
& \hspace{4mm} + \left( (\mu +\delta) + \dfrac{a_3}{\vartheta (\nu + \eta + \delta)} \right) \left( \mu (\nu + \eta + 2\delta) + \dfrac{a_3}{\vartheta} + \dfrac{a_3}{\nu + \eta + \delta} \right) > 0.
\end{align*}
Finally, the conditions \eqref{A-RHsys1} are satisfied and we conclude that $\Eb_1^{*}$ is LAS whenever $\N_F > 1.$

\textbf{Proof of Proposition \ref{prop3}.} The Jacobian related to the system \eqref{sys-F2} with $\Phib_2$ defined by \eqref{F2} is given by

\[ J(\Xb)= \begin{pmatrix}
 -\mu \!-\! r  \rho \sigma U e^{-\sigma(A+M+U)}  & - r \rho \sigma U e^{-\sigma (A+M+U)} & r \rho (1 \!-\! \sigma  U) e^{-\sigma  (A+M+U)} \\
& & \\
 -\gamma  \nu \!-\! (1 \!-\! r) \rho \sigma U e^{-\sigma (A+M+U)}  & -\delta \!-\! (1-r) \rho \sigma U e^{\sigma-(A+M+U)} & \eta \!+\!  (1 \!-\! r) \rho (1 \!-\! \sigma U) e^{-\sigma  (A+M+U)} \\
& & \\
 \gamma  \nu  & 0 & -(\eta \!+\! \delta ) \end{pmatrix} \]
where $\Xb=(M,A,U)$. Thus, it is easy to show that

\[ J(\Eb_0)=\begin{pmatrix}
 -\mu  & 0 & r \rho  \\
 -\gamma  \nu  & -\delta  & \eta +(1-r) \rho  \\
 \gamma  \nu  & 0 & -\delta -\eta  \\
\end{pmatrix}, \]
and to show that the characteristic polynomial is given by

\[ P_2^0(\lambda)=-(\delta + \lambda ) \Big[ \lambda^2 + (\delta + \eta + \mu) \lambda  + (\delta +\eta) \mu  -\gamma  \nu  \rho  r \Big] = -(\delta + \lambda ) \Big[ \lambda^2 + (\delta + \eta + \mu) \lambda  + (\delta +\eta) \mu \big( 1 - \N_M \big) \Big] \]

When $\N_M<1$,  this polynomial has three roots with negative real parts meaning that $\Eb_0$ is LAS. Alternatively, $p_2(\lambda)$ has one root with positive real part when $\N_M > 1$ meaning that $\Eb_0$ is a saddle point (not a repeller). In effect, a trajectory engendered by the initial condition $M(0) >0, A(0)=0, U(0)=0$ converges to $\Eb_0$ even if $\N_M >1$.

Let us now show that $E_2^*$ is LAS  when $\N_M>1$. We recall here that $U_2^*=\dfrac{\gamma \nu}{\eta + \delta} M_2^*,$ and thus we have

\begin{equation}
\label{A-JE2}
J \big(\Eb_2^* \big)=\begin{pmatrix}
-\mu \big(1 + \sigma M^*_2 \big)                     & -\mu \sigma  M^*_2             & \dfrac{r\rho}{\N_M} -\mu \sigma M_2^*\\
& & \\
- \left( \gamma  \nu + \dfrac{(1-r)  }{r} \mu \sigma  M^*_2 \right) & - \left( \delta + \dfrac{(1-r)}{r} \mu \sigma M^*_2 \right) & \eta + \dfrac{(1-r)\rho}{\N_M} - \dfrac{1-r}{r} \mu \sigma M_2^* \\
& & \\
\gamma \nu & 0 & - (\eta+\delta )
\end{pmatrix}
\end{equation}

The characteristic polynomial of $J \big(\Eb_2^* \big)$ has the form

\[ P_2^* (\lambda) =  \lambda^3 + b_1 \lambda^2 +b_2 \lambda + b_3 ,  \]
where $b_1 = - \text{\texttt{Tr} } J \big(\Eb_2^{*} \big), b_3 = - \det J \big(\Eb_2^{*} \big),$ and

\begin{equation}
 \label{A-b2}
b_2 = \det \begin{vmatrix} J_{11} & J_{12} \\ J_{21} & J_{22} \end{vmatrix} + \det \begin{vmatrix} J_{11} & J_{13} \\ J_{31} & J_{33} \end{vmatrix} + \det \begin{vmatrix} J_{22} & J_{23} \\ J_{32} & J_{33} \end{vmatrix} := \Delta_1^{(2)} + \Delta_2^{(2)} + \Delta_3^{(2)}
\end{equation}
with $J_{ij}, i,j=1,2,3$ denoting the elements of \eqref{A-JE2}. According to Routh-Hurwitz criterion (see, e.g., \cite{Murray2002}), all roots of $P_2^*$ have negative real parts if and only if the following conditions are satisfied:

\begin{equation}
\label{A-RHsys2}
b_1 >0, \quad b_3 > 0,   \quad \text{and} \quad b_1 b_2 - b_3 >0.
\end{equation}
Let us now check these conditions. First we note that

\[ \text{\texttt{Tr} } J \big(\Eb_2^{*} \big) = - 2 \delta - \eta - \mu  \big( 1 + \sigma M^*_2 \big) - \frac{(1-r)}{r} \mu  \sigma M^*_2 < 0 , \]
and thus we have

\[ b_1 = - \text{\texttt{Tr} } J \big(\Eb_2^{*} \big) = \mu + \eta + 2\delta + \dfrac{\mu}{r} \sigma M^*_2 >0.
\]

To compute $\det J\big( \Eb_2^{*} \big)$, we observe that

\[  \det J \big( \Eb_2^{*} \big) = \gamma \nu \Delta_{31}^{(2)} - (\eta+\delta ) \Delta_{33}^{(2)}, \]
where $\Delta_{31}^{(2)}, \Delta_{33}^{(2)}$ are minors of $J \big( \Eb_2^{*} \big)$ obtained by elimination of the third row and either first or third column from \eqref{A-JE2}. Effectively, we have

\begin{align*}
\Delta_{31}^{(2)} & = - \mu \sigma M_2^* \left( \eta + \dfrac{(1-r)\rho}{\N_M} - \dfrac{1-r}{r} \mu \sigma M_2^* \right) + \left( \dfrac{r\rho}{\N_M} -\mu \sigma M_2^* \right) \left( \delta + \dfrac{(1-r)}{r} \mu \sigma M^*_2 \right) \\
&  = - \eta \mu \sigma M_2^* - \delta \mu \sigma M_2^* + \delta \dfrac{r \rho}{\N_M} = - (\eta + \delta) \mu \sigma M_2^* + \dfrac{\mu \delta (\eta + \delta)}{\gamma \nu} = \mu (\eta + \delta) \left[ \dfrac{\delta}{\gamma \nu} - \sigma M_2^* \right], \\
\Delta_{33}^{(2)} & = \mu \big(1 + \sigma M^*_2 \big) \left( \delta + \dfrac{(1-r)}{r} \mu \sigma M^*_2 \right) - \mu \sigma  M^*_2 \left( \gamma  \nu + \dfrac{(1-r)  }{r} \mu \sigma  M^*_2 \right) \\
& = \delta \mu + \left( \delta + \mu \dfrac{(1-r)}{r} \right) \mu \sigma M_2^*   - \gamma \nu \mu \sigma M_2^* = \delta \mu + \dfrac{\vartheta \mu }{r} \sigma M_2^* - \gamma \nu \mu \sigma M_2^*,
\end{align*}
and therefore

\begin{align*}
\det J \big( \Eb^*_2 \big)  & =   \gamma \nu \mu (\eta + \delta) \left[ \dfrac{\delta}{\gamma \nu} - \sigma M_2^* \right] - (\eta + \delta) \left[ \delta \mu + \dfrac{\vartheta \mu }{r} \sigma M_2^* - \gamma \nu \mu \sigma M_2^* \right] \\
 & = (\eta + \delta) \left[ \delta \mu - \gamma \nu \mu \sigma M_2^* - \delta \mu - \dfrac{\vartheta }{r} \mu \sigma M_2^* + \gamma \nu \mu \sigma M_2^*\right] = - (\eta + \delta) \dfrac{\vartheta}{r}  \mu \sigma M_2^* < 0\\
\end{align*}

Thus we have

\begin{equation}
\label{A-b3}
b_3 = - \det J \big( \Eb^*_2 \big) = (\eta + \delta) \dfrac{\vartheta}{r}  \mu \sigma M_2^* > 0
\end{equation}

To compute the coefficient $b_2$, we evaluate $\Delta_i^{(2)}, i=1,2,3$ that appear in \eqref{A-b2}:

\begin{align*}
\Delta_1^{(2)} & = \Delta_{33}^{(2)} =  \delta \mu + \dfrac{\vartheta \mu }{r} \sigma M_2^* - \gamma \nu \mu \sigma M_2^*, \\
\Delta_2^{(2)} & = (\eta + \delta) \mu \big(1 + \sigma M^*_2 \big) - \gamma  \nu \left( \dfrac{r\rho}{\N_M} - \mu \sigma M^*_2 \right) = (\eta + \delta) \mu \big(1 + \sigma M^*_2 \big) - \mu (\eta + \delta) + \gamma  \nu \mu \sigma  M^*_2   \\
 & = (\eta + \delta) \mu \sigma M^*_2 + \gamma  \nu \mu \sigma  M^*_2, \\
\Delta_3^{(2)} & = (\eta + \delta) \left( \delta + \dfrac{(1-r)}{r} \mu \sigma M^*_2 \right).
\end{align*}
Then using \eqref{A-b2}, we obtain

\begin{align*}
b_2 & = \delta \mu + \dfrac{\vartheta \mu }{r} \sigma M_2^* + (\eta + \delta) \mu \sigma M^*_2 + \delta (\eta + \delta) + (\eta + \delta) \dfrac{(1-r)}{r} \mu \sigma M^*_2 \\
& = \delta ( \mu  + \eta + \delta) +  (\eta + \delta) \dfrac{1}{r} \mu \sigma M^*_2 + \dfrac{\vartheta}{r} \mu \sigma M^*_2 = \delta ( \mu  + \eta + \delta) + \dfrac{\eta + \delta + \vartheta}{r} \mu \sigma M^*_2 >0.
\end{align*}

Our goal now is to show that $b_1 b_2 - b_3 > 0.$ Before proceeding, we rewrite $b_1$ and $b_2$ in terms of $b_3$ using the formula \eqref{A-b3}:

\[ b_1 = (\mu + \delta) + (\eta + \delta) + \dfrac{b_3}{\vartheta (\eta + \delta)}, \qquad b_2 = \delta (\mu  + \eta + \delta) + \dfrac{b_3}{\vartheta} + \dfrac{b_3}{\eta + \delta}, \]
so that

\begin{align*}
b_1 b_2 - b_3 & = \left( (\mu + \delta)  + \dfrac{b_3}{\vartheta (\eta + \delta)} + (\eta + \delta) \right) \left( \delta (\mu  + \eta + \delta) + \dfrac{b_3}{\vartheta} + \dfrac{b_3}{\eta + \delta} \right) - b_3 \\
& = (\eta + \delta) \left( \delta (\mu  + \eta + \delta) + \dfrac{b_3}{\vartheta} \right) + b_3 + \left( (\mu + \delta)  + \dfrac{b_3}{\vartheta (\eta + \delta)} \right) \left( \delta (\mu  + \eta + \delta) + \dfrac{b_3}{\vartheta} + \dfrac{b_3}{\eta + \delta} \right)- b_3 \\
& = (\eta + \delta) \left( \delta (\mu  + \eta + \delta) + \dfrac{b_3}{\vartheta} \right) + \left( (\mu + \delta)  + \dfrac{b_3}{\vartheta (\eta + \delta)} \right) \left( \delta (\mu  + \eta + \delta) + \dfrac{b_3}{\vartheta} + \dfrac{b_3}{\eta + \delta} \right) >0.
\end{align*}
Finally, the conditions \eqref{A-RHsys2} are satisfied and  we conclude that  $\Eb_2^{*}$ is LAS whenever $\N_M > 1.$
\end{appendix}

\begin{appendix}
\appendix
\renewcommand{\theequation}{B-\arabic{equation}}
  \setcounter{equation}{0}  
\renewcommand{\thefigure}{B.\arabic{figure}}
  \setcounter{figure}{0}  
  \renewcommand{\thetable}{B.\arabic{table}}
   \setcounter{table}{0}  

\section*{Appendix B: proof of Proposition \ref{existence_equilibrium_scarse_pheromones}}
\addcontentsline{toc}{section}{Appendix B: proof of Proposition \ref{existence_equilibrium_scarse_pheromones}}
\label{appendixB}

We have to solve the following algebraic system

\begin{subequations}
\label{syst_to_solve}
   \begin{align}[left = \empheqlbrace\,]
    \label{syst_to_solve1}
    r\rho Ue^{-\sigma(M+A+U)}-\alpha\dfrac{A_{p}}{A+A_{p}}M-\mu M &=0, \\[1mm]
     \label{syst_to_solve2}
(1-r)\rho Ue^{-\sigma(M+A+U)}-\gamma\nu\dfrac{A}{A+A_{p}}M+\eta U-\delta A &=0,\\[1mm]
 \label{syst_to_solve3}
\gamma\nu\dfrac{A}{A+A_{p}}M-\eta U-\delta U &=0.
    \end{align}
\end{subequations}
From Eqs. \eqref{syst_to_solve1} and \eqref{syst_to_solve3}, we obtain first

\begin{equation}
\label{sys_12}
\rho Ue^{-\sigma(M+A+U)}=\dfrac{1}{r}\left(\alpha\dfrac{A_{p}}{A+A_{p}}+\mu\right)M, \qquad
U=\dfrac{\gamma\nu}{\eta+\delta} \: \dfrac{A}{A+A_{p}} M
\end{equation}
and then using \eqref{syst_to_solve2}, we arrive to

\[ \left( \dfrac{1-r}{r} \left(\alpha \dfrac{A_{p}}{A+A_{p}} + \mu \right) - \gamma \nu \dfrac{A}{A+A_{p} } + \eta \dfrac{\gamma \nu}{\eta + \delta} \dfrac{A}{A+A_{p}} \right) M = \delta A. \]
The above expression can also be written as

\[ \left( \left( 1-r \right) \left( \alpha + \mu \right) A_{p }+ \left( \left( 1-r \right) \mu - \dfrac{\delta \gamma \nu r}{\eta + \delta} \right) A \right) M = r \delta \big( A + A_{p} \big) A. \]

Further, using the quantity $\theta_M$, defined by \eqref{theta_M}, we arrive to

\[ \left( \left( 1-r \right) \left( \alpha + \mu \right) A_{p} + \dfrac{\delta \gamma \nu r}{\eta + \delta} \big( \theta_M -1 \big) A \right) M = r \delta \big( A + A_{p}\big) A. \]
Here, when $\theta_M>1$, we have

\[ M = \dfrac{r \delta \big( A+A_{p} \big)}{\left( 1 - r \right) \left( \alpha + \mu \right) A_{p} + \dfrac{\gamma r \nu \delta }{\eta + \delta} \big( \theta_M - 1 \big) A} A, \]
or

\[ M = \dfrac{\eta + \delta}{\gamma \nu} \: \dfrac{ A + A_{p} }{\dfrac{( \eta + \delta) ( 1 - r ) ( \alpha + \mu )}{\gamma r \nu \delta} + \big( \theta_M - 1 \big) A} A \]
and

\[ U = \dfrac{\gamma \nu}{\eta +\delta} \: \dfrac{A}{A + A_{p}} M = \dfrac{1}{\dfrac{ \left(\eta + \delta\right) \left(1 - r \right) \left( \alpha + \mu \right)}{\gamma r \nu \delta} A_{p} + \big( \theta_{M} - 1 \big) A}A^{2}.
\]
Thus, $M$ and $U$ are now expressed in terms of $A$ and the external parameters $\alpha, A_p$. By plugging the above expressions for $M$ and $U$ into the left-hand side formula of \eqref{sys_12}, we obtain the following relationship

\[ A e^{-\sigma (M+A+U)} = \dfrac{1}{\rho r}\Big( (\alpha + \mu)A_{p}+\mu A \Big) \dfrac{\eta + \delta}{\gamma \nu}=\dfrac{1}{\mu \N_M} \Big( (\alpha + \mu)A_{p} +\mu A \Big). \]

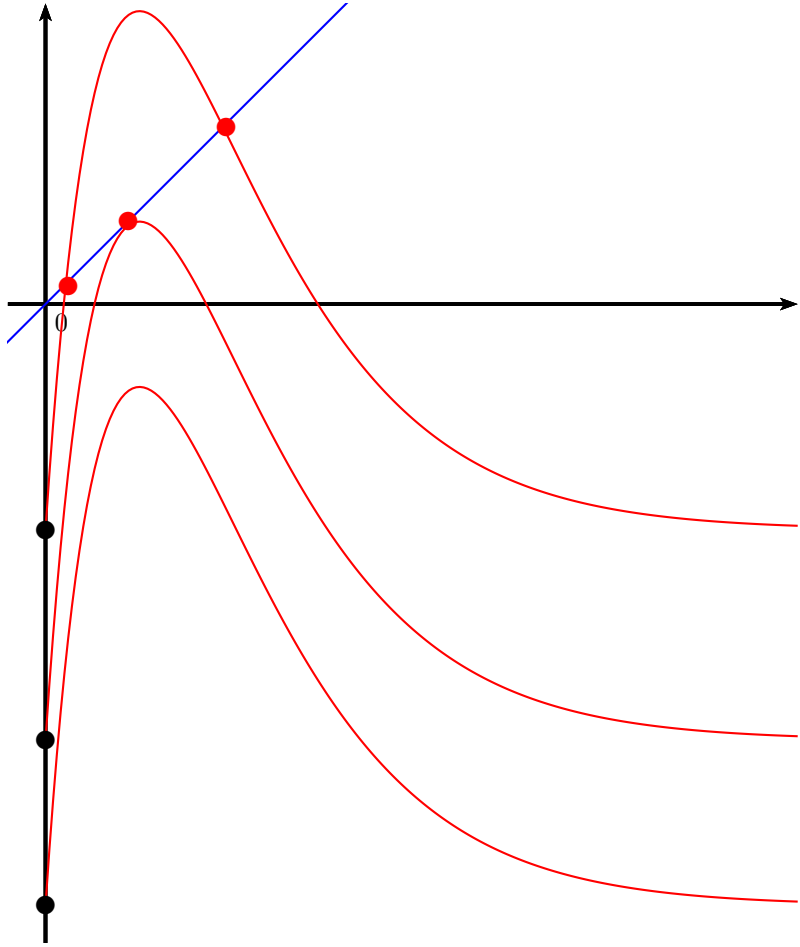
\begin{figure}[t]
\centering

\newrgbcolor{wwwwww}{0.4 0.4 0.4}
\begin{pspicture*}(-0.5,-8.5)(10,4)
\psset{xunit=1.0cm,yunit=1.0cm,algebraic=true,dotstyle=*,dotsize=7pt 0,linewidth=0.8pt,arrowsize=2pt 2,arrowinset=0.25}
\psaxes[labelFontSize=\scriptstyle,xAxis=true,yAxis=true,linewidth=1.4pt,Dx=20,Dy=20,ticksize=-2pt 0,subticks=2]{->}(0,0)(-0.5,-8.5)(10,4)  
\uput[dr](0,0){$0$}
\psplot[linecolor=red,plotpoints=200]{0}{10}{15*x*EXP(-0.8*x)-8}
\psplot[linecolor=red,plotpoints=200]{0}{10}{15*x*EXP(-0.8*x)-5.8}
\psplot[linecolor=red,plotpoints=200]{0}{10}{15*x*EXP(-0.8*x)-3}
\psplot[linecolor=blue,plotpoints=200]{-1}{11}{x}
\uput[u](10.75,0.25){$A$}
\psdot(0,-3)
\psdot(0,-8)
\psdot(0,-5.8)
\psdot[linecolor=red](1.1,1.1)
\psdot[linecolor=red](0.3,0.24)
\psdot[linecolor=red](2.4,2.35)
\end{pspicture*}
\caption{Possible intersections between $\varphi(A;A_p)$ (red color) and $A$ (blue color) for different values of $A_p$ }
\label{fig_case}
\end{figure}

Finally, replacing $M$ and $U$ in the exponential term leads to the following equation to solve

\begin{equation}
\label{eqA}
\N_M A f(A;A_p)-\dfrac{\alpha+\mu}{\mu}A_{p}=A,
\end{equation}
where

\[ f(A;A_p) := \exp \left(- \sigma \left( 1 + \dfrac{\delta r A_{p} + \delta r \left( \dfrac{\gamma \nu}{\eta + \delta} + 1 \right) A}{\left( 1 - r \right) \left( \alpha + \mu \right) A_{p} + \dfrac{\delta r \gamma \nu}{ \eta + \delta } \big(\theta_{M} - 1 \big) A} \right) A \right),
\]
which is a function of $A$ also depending on the external parameter $A_p \geq 0$. Notably, when $A_p=0$, we recover $f(A ; 0) = \dfrac{1}{\N_M}$. For the fixed values of $\alpha$ and $A_p$, let us denote the left-hand side of the equation \eqref{eqA} by the function

\[ \varphi(A;A_p) := \N_M A f(A ; A_p ) - \dfrac{\alpha + \mu}{\mu} A_{p} \]
that fulfills the condition $\varphi(0;A_p) <0$ whenever $A_p >0$. When $A>0$, function $\varphi(A;A_p)$ increases first and then decreases. Therefore, depending on the value $A_p >0$, there may exist two, one, or no solutions to equation \eqref{eqA} whose right-hand side is a straight line, see Figure \ref{fig_case}. As shown in Figure \ref{fig_case}, equation \eqref{eqA} has two solutions when $A_p$ is relatively small. Then, by gradually increasing the value of $A_p$, one may get only one solution of \eqref{eqA}. In such a case, the corresponding value of $A_p$ will render the threshold value $A_p > A_p^{crit}$. Namely, equation \eqref{eqA} has no solution when $A_p > A_p^{crit}$ and no positive equilibria of the system \eqref{syscon_scarse_open_loop} can exist for $A_p > A^{crit}_p$.

Thus, to identify the threshold value $A_p^{crit}$ together with the underlying unique solution to equation \eqref{eqA}, one must resolve the system of two equations:

\begin{equation}
\label{eqA-Ap}
\varphi(A;A_p)=A, \qquad \dfrac{\partial \varphi (A;A_p)}{\partial A} =1.
\end{equation}
Even though no analytical formula can be obtained for $A_p^{crit}$, its underlying value can be adequately approximated by solving the system \eqref{eqA-Ap} numerically.

\end{appendix}

\begin{appendix}
\appendix
\renewcommand{\theequation}{C-\arabic{equation}}
  \setcounter{equation}{0}  
\renewcommand{\thefigure}{C.\arabic{figure}}
  \setcounter{figure}{0}  
  \renewcommand{\thetable}{C.\arabic{table}}
   \setcounter{table}{0}  

\section*{Appendix C: proof of Theorem \ref{theo_GAS_scarse_aux} }
\addcontentsline{toc}{section}{Appendix C: proof of Theorem \ref{theo_GAS_scarse_aux}}
\label{appendixC}

The formal proofs of items (a) and (b) are similar to previous explanations, so we leave them to the readers. Let us focus on item (c). Setting the first and third components of \eqref{auxiliary_scarse_open_loop} equal to zero leads to

\[ r \rho U= \left( \alpha \dfrac{A_{p}}{A + A_{p}} + \mu \right) M e^{\sigma M}, \]
and, assuming $M>0$,

\[ r \rho \gamma\nu A - ( \eta + \delta) \Big( (\alpha + \mu  )A_{p}+\mu A \Big) e^{\sigma M} = 0, \]
that is

\[ \big( \N_{M} - e^{\sigma M} \big) A = \dfrac{\alpha + \mu}{\mu} A_{p} e^{\sigma M} \quad \Rightarrow \quad A = \dfrac{(\alpha + \mu)}{\mu \big( \N_M - e^{\sigma M} \big)} A_p e^{\sigma M}. \]

Setting the second component of \eqref{auxiliary_scarse_open_loop} to zero and replacing $A$ and $U$ provides

\begin{align*}
r \rho \delta \dfrac{\alpha + \mu}{\mu} A_{p} &= \big(\N_{M}-e^{\sigma M}\big) \Big( (1-r)\rho + \eta \Big) \left( \alpha \dfrac{A_{p}}{A+A_{p}} + \mu \right) M \\
&= \big( \N_{M} - e^{\sigma M} \big) \Big( (1-r)\rho+\eta \Big) \left( \alpha \dfrac{A_{p}}{\dfrac{\alpha+\mu}{\mu \big( \N_{M} - e^{\sigma M} \big)} A_p e^{\sigma M} + A_{p}} + \mu \right) M, \\[2mm]
&= \big( \N_{M} - e^{\sigma M} \big) \Big( (1-r)\rho+\eta \Big) \mu \dfrac{(\alpha + \mu) \N_M}{\alpha e^{\sigma M} + \mu N_M} M.
\end{align*}
From the above relationships, we have the following equation

\begin{equation}
    \label{auxiliary_system_zero}
\Big( (1-r)\rho+\eta \Big) \mu \N_M M \big( \N_{M} - e^{\sigma M} \big) = r \rho \delta A_{p}\left( \N_{M} + \dfrac{\alpha}{\mu}e^{\sigma M} \right).
\end{equation}
which can be viewed as $g_1(M) = g_2(M; A_p,\alpha),$ where

\[ g(M) := \Big( (1-r)\rho+\eta \Big) \mu \N_M M \big( \N_{M} - e^{\sigma M} \big) \quad \text{and} \quad h(M;A_p,\alpha) :=   r \rho \delta A_{p}\left( \N_{M} + \dfrac{\alpha}{\mu}e^{\sigma M} \right). \]

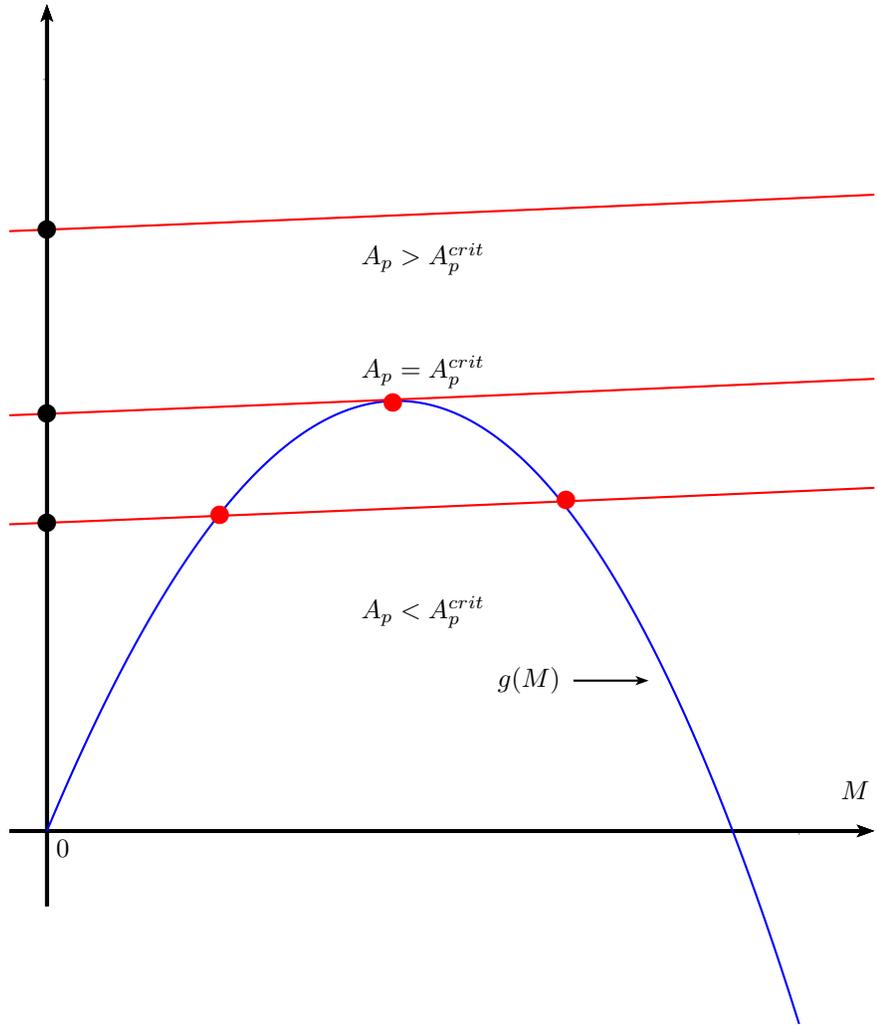
\begin{figure}[t]
\centering

\newrgbcolor{wwwwww}{0.4 0.4 0.4}
\begin{pspicture*}(-0.5,-4)(11,11)
\psset{xunit=1.0cm,yunit=1.0cm,algebraic=true,dotstyle=*,dotsize=7pt 0,linewidth=0.8pt,arrowsize=2pt 2,arrowinset=0.25}
\psaxes[labelFontSize=\scriptstyle,xAxis=true,yAxis=true,linewidth=1.4pt,Dx=20,Dy=20,ticksize=-2pt 0,subticks=2]{->}(0,0)(-0.5,-1)(11,11)  
\uput[dr](0,0){$0$}
\psplot[linecolor=blue,plotpoints=200]{0}{10}{12*x*(1.2-EXP(0.02*x))}
\psplot[linecolor=red,plotpoints=200]{-1}{11}{1.55+4*EXP(0.01*x)}
\psplot[linecolor=red,plotpoints=200]{-1}{11}{4.+4*EXP(0.01*x)}
\psplot[linecolor=red,plotpoints=200]{-1}{11}{0.1+4*EXP(0.01*x)}
\psdot(0,8)
\psdot(0,5.55)
\psdot(0,4.1)
\uput[u](10.75,0.25){$M$}
\psdot[linecolor=red](4.6,5.7)
\psdot[linecolor=red](6.9,4.4)
\psdot[linecolor=red](2.3,4.2)
\uput[d](5,3.3){$A_{p}<A_{p}^{crit}$}
\uput[d](5,6.5){$A_{p}=A_{p}^{crit}$}
\uput[d](5,8){$A_{p}>A_{p}^{crit}$}
\uput[l](7,2){$g(M)$}
\psline{->}(7,2)(8,2)
\end{pspicture*}

\caption{Possible intersections between \{$g(M)$ (blue color) and $h(M;A_p,\alpha)$ (red color) for different values of $A_p$ \label{fig:sys-reduit}}
\end{figure}

Here, the function $g(M)$ is increasing for small values of $M$ and then decreasing as $M$ becomes larger, so there exists $\widehat{M} >0$ where $g_1$ attains a maximum $g \big( \widehat{M} \big)$. On the other hand, the function $h(M; A_p,\alpha)$ is increasing for all $A_p >0$ and $\alpha >0$ (or constant when $\alpha =0$). Thus, for any fixed $\alpha \in [0,1]$ there exists a quantity $\tilde{A}_p^{crit} >0$ such that Eq. \eqref{auxiliary_system_zero} has: (a) no positive roots when $A_p > \tilde{A}_p^{crit}$; (b) two positive roots $M_{1,p}^*$ and $M_{2,p}^*$ when $A_p < \tilde{A}_p^{crit}$. Notably, when $A_p = \tilde{A}_p^{crit}$, the two positive roots collide ($M_{1,p}^*=M_{2,p}^*$) meaning that the auxiliary system \eqref{auxiliary_scarse_open_loop}  undergoes a pitchfork bifurcation. Figure \ref{fig:sys-reduit} illustrates possible intersections between $g(M)$ (blue-colored curve) and $g(M; A_p,\alpha)$ (red-colored curves) for different values of $A_p$. In the sequel, we address items (a)  and (b) mentioned above in more detail.

\begin{itemize}
    \item[(a)]
    Assume $A_{p}>\tilde{A}_{p}^{crit}$ and let $q \in \mathbb{R}_{+}$, with $q \geq q^{*} := \dfrac{\mu}{r \rho} \dfrac{\N_{M}}{\sigma} \ln \N_{M}$. We denote the right-hand side of the auxiliary system \eqref{auxiliary_scarse_open_loop} by
    
\[ \Hb(\Xb) :=\begin{pmatrix}
r \rho U e^{-\sigma M} - \alpha \dfrac{A_{p}}{A_{p} + A}M - \mu M \\[4mm]
(1-r)\rho U + \eta U - \delta A \\[3mm]
\nu \dfrac{\gamma M}{A_{p} + A} A - \eta U - \delta U
\end{pmatrix}, \]
where $\Xb=(M,A,U)$, and define

\[ \Xb_{q} := \begin{pmatrix}
\dfrac{r \rho}{\mu} \dfrac{1}{N_{M}} q \\[4mm]
2 \dfrac{(1-r) \rho + \eta}{\delta} q \\[3mm]
q
\end{pmatrix}. \]
It is not difficult to check that  $\Hb \big( \Eb_{0} \big)= \Eb_{0}$ and $\Hb \big( \Xb_{q} \big) \leq \Eb_{0}$. Then, according to Anguelov \textit{et al.} (see Theorem 7 in \cite{Anguelov2017}), we deduce that $\Eb_{0}$ is globally asymptotically stable (GAS) on $\big[ \Eb_{0}, \Xb_{q} \big]$. Therefore, $\Eb_{0}$ is GAS on $\Omega$ and also on $\mathbb{R}_{+}^{3}$ because $\Omega$ is an absorbing set (see Proposition \ref{prop1}).

\item[(b)]
Assume $A_{p}<\tilde{A}_{p}^{crit}$. Let $\tilde{\Eb}_{1,p}^*$ and $\tilde{\Eb}_{2,p}^*$ be two equilibria such that $\tilde{\Eb}_{i,p}^*=\big( M_{i,p}^*,A_{i,p}^*,U_{i,p}^*\big)$, $i=1,2$, where $M_{1,p}^*$ and $M_{2,p}^*$ are positive roots of \eqref{auxiliary_system_zero} that fulfill $M_{1,p}^* < M_{2,p}^*$ (in an ``element-by-element'' sense). Since $A$ and $U$ are increasing functions of $M$, we deduce that $\Eb_0 < \tilde{\Eb}_{1,p}^* < \tilde{\Eb}_{2,p}^*$. As it holds that $\Hb \big(\Xb_{q} \big) \leq \Eb_{0}= \Hb \big( \tilde{\Eb}_{2,p}^* \big)$, by applying again Theorem 7 from \cite{Anguelov2017}, we conclude that $\tilde{\Eb}_{2,p}^*$ is GAS on $\big[ \tilde{\Eb}_{2,p}^*,\Xb_q \big]$, that is, $\tilde{\Eb}_{2,p}^*$ is GAS on $\bigcup \limits_{q \geq q^*} \big[ \tilde{\Eb}_{2,p}^*,\Xb_q \big] = \Big\{ \Xb \in \mathbf{R}_+^3: \ \Xb \geq \tilde{\Eb}_{2,p}^* \Big\}$.

Finally, we deduce a similar result for $\Eb_0 <\tilde{\Eb}_{1,p}^*$, using Theorem 8 from \cite{Anguelov2017}, namely, all solutions initiated in $ \Big\{ \Xb \in \mathbf{R}_+^3: \ \Xb <  \tilde{\Eb}_{1,p}^* \Big\}$ converge to $\Eb_0$.
\end{itemize}

\end{appendix}

\end{document}